\documentclass[12pt,a4paper]{amsart}
\usepackage{amsfonts}
\usepackage{amsthm}
\usepackage{amsmath}
\usepackage{amssymb}
\usepackage{amscd}
\usepackage[latin2]{inputenc}
\usepackage{t1enc}
\usepackage[mathscr]{eucal}
\usepackage{indentfirst}
\usepackage{graphicx}
\usepackage{graphics}
\usepackage{pict2e}
\usepackage{epic}
\numberwithin{equation}{section}
\usepackage[margin=2.9cm]{geometry}
\usepackage{epstopdf} 

\usepackage{multirow} 
\usepackage{array}

\usepackage{soul}
\usepackage[dvipsnames]{xcolor}

\theoremstyle{definition}
\newtheorem{theorem}{Theorem}[section]
\newtheorem{lemma}[theorem]{Lemma}
\newtheorem{corollary}[theorem]{Corollary}
\newtheorem{proposition}[theorem]{Proposition}
\newtheorem{definition}[theorem]{Definition}
\newtheorem{example}[theorem]{Example}

\newtheorem{Rem}[theorem]{Remark}

\newcommand{\im}{\operatorname{im}}

\begin{document}

\title[Geometric bijections between subgraphs and orientations of a graph]{Geometric bijections between spanning subgraphs and orientations of a graph}

\author[Changxin Ding]{Changxin Ding}

\address{School of Mathematics, Georgia Institute of Technology \\ Atlanta, Georgia 30332-0160, USA } 

\email{cding66@gatech.edu}

\begin{abstract}
Let $G$ be a connected finite graph. Backman, Baker, and Yuen have constructed a family of explicit and easy-to-describe bijections between spanning trees of $G$ and $(\sigma,\sigma^*)$-compatible orientations, where the $(\sigma,\sigma^*)$-compatible orientations are the representatives of equivalence classes of orientations up to cycle-cocycle reversal which are determined by a cycle signature $\sigma$ and a cocycle signature $\sigma^*$. Their bijections are geometric since the construction comes from zonotopal subdivisions.

In this paper, we extend the geometric bijections to subgraph-orientation correspondences. Moreover, we extend the geometric constructions accordingly. Our proofs are combinatorial, which do not make use of the zonotopes. We also provide geometric proofs 
for partial results, which make use of zonotopal tiling, relate to Backman, Baker, and Yuen's method, and motivate our combinatorial constructions. Finally, we explain that the main results hold for \emph{regular matroids}.

\end{abstract}

\maketitle


\section{Introduction}\label{intro}
\subsection{Introduction to the main combinatorial results}
Let $G$ be a connected finite graph with nonempty edge set $E$, where \emph{loops} and \emph{multiple edges} are allowed.
This paper examines correspondence between spanning subgraphs and orientations of $G$. Obviously, the number of spanning subgraphs of $G$ equals the number of orientations of $G$. More interestingly, some types of spanning subgraphs are equinumerous to some types of equivalence classes of orientations, which can be counted by the \emph{Tutte polynomial} $T_G(x,y)$ of $G$; see Table~\ref{table1}. The Tutte polynomial counts the types of spanning subgraphs in the first column by definition. The equivalence classes of orientations in the second column were introduced and enumerated by Gioan \cite{G1}. It is therefore natural to look for bijections between them.

\begin{table}[h!]
\centering
\begin{tabular}{ |m{5.5cm}|m{6.5cm}|m{1.8cm}| }
\hline
types of spanning subgraphs & equivalence classes of orientations & cardinality\\
\hline
spanning subgraphs  &   orientations   & $T_G(2,2)$\\
spanning forests  &  equivalence classes of orientations up to cycle reversal   & $T_G(2,1)$  \\
connected spanning subgraphs &  equivalence classes of orientations up to cocycle reversal & $T_G(1,2)$\\
spanning trees   & equivalence classes of orientations up to cycle-cocycle reversal & $T_G(1,1)$\\
 
 \hline
\end{tabular}
\caption{Some spanning subgraphs and orientations counted by the Tutte polynomial. By spanning forests we mean acyclic subsets of edges.}
\label{table1}
\end{table}

\begin{figure}[h]
            \centering
            \includegraphics[scale=1]{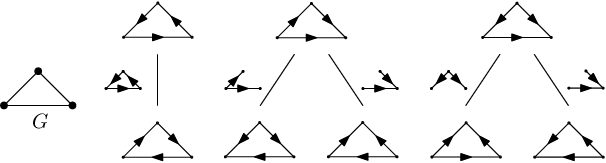}
            \caption{The triangle graph $G$ and its eight orientations. They are partitioned into the three cycle-cocycle reversal classes. From the top three orientations, one may get the others by reversing certain directed cycles or cocycles, which are indicated beside the line connecting two orientations.}
            \label{F-1}
\end{figure}

In~\cite{BBY}, Backman, Baker, and Yuen construct a family of bijections for the last row of Table~\ref{table1}. We denote by $\mathcal{T}(G)$ the set of spanning trees of $G$. We call the equivalence classes of orientations up to cycle-cocycle reversal the \emph{cycle-cocycle reversal (equivalence) classes}, the set of which is denoted by $\mathcal{G}(G)$. Two orientations are in the same class if and only if one can obtain one orientation from the other by reversing some directed cycles and cocycles; see Section~\ref{combinatorial} for details and see Figure~\ref{F-1} for an example.

In \cite{BBY}, the authors actually establish a bijection from $\mathcal{T}(G)$ to the set of \emph{$(\sigma,\sigma^*)$-compatible orientations} of $G$, where the latter set is a representative set of $\mathcal{G}(G)$ determined by a pair of \emph{acyclic cycle signature} $\sigma$ and \emph{acyclic cocycle signature} $\sigma^*$.  In this paper, we call this bijection the \emph{BBY bijection} and denote it by $\text{BBY}_{\sigma,\sigma^*}$. For a fixed graph, we can choose different pairs of acyclic signatures $\sigma$ and $\sigma^*$ and get different BBY bijections. The BBY bijections are geometric since their construction comes from \emph{zonotopal subdivisions}. 

Let us mention that their work is motivated by finding bijections between $\mathcal{T}(G)$ and the \emph{Jacobian group} (also known as the \emph{critical group}, \emph{sandpile group}, or \emph{Picard group}) $\text{Jac}(G)$ of $G$. There is a canonical simply transitive group action of $\text{Jac}(G)$ on $\mathcal{G}(G)$ \cite{B}, and hence the BBY bijection induces a bijection between $\mathcal{T}(G)$ and $\text{Jac}(G)$.

Our main work is to extend the BBY bijection to a subgraph-orientation correspondence. To present their bijections and ours, we now give the necessary definitions and properties.

\begin{figure}[h]
            \centering
            \includegraphics[scale=1]{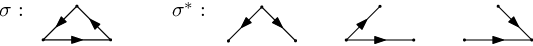}
            \caption{An acyclic cycle signature $\sigma$ and an acyclic cocycle signature $\sigma^*$ of the triangle graph.}
            \label{F-2}
\end{figure}

A \emph{cycle signature} $\sigma$ is the choice of a direction for each cycle of the graph $G$. For each cycle $C$, we denote by $\sigma(C)$ the directed cycle we choose for $C$. By abusing notation, we also view $\sigma$ as the set of the directed cycles we choose: $\{\sigma(C):C\text{ is a cycle}\}$. By fixing a reference orientation of $G$, we can identify directed cycles with $\{0,\pm 1\}$-vectors in $\mathbb{R}^E$. The cycle signature $\sigma$ is said to be \emph{acyclic} if whenever $a_C$ are nonnegative reals with $\sum_C a_C\sigma(C)=0$ in $\mathbb{R}^E$ we have $a_C=0$ for all $C$, where the sum is over all cycles of $G$. An orientation is said to be \emph{$\sigma$-compatible} if any directed cycle in the orientation is in $\sigma$.  An \emph{acyclic cocycle signature} $\sigma^*$ and a \emph{$\sigma^*$-compatible orientation} are defined similarly but for directed cocycles instead of directed cycles.  An orientation is said to be  \emph{($\sigma,\sigma^*$)-compatible} if it is both $\sigma$-compatible and $\sigma^*$-compatible.
For any orientation $\overrightarrow{O}$ of $G$, there exists a unique $(\sigma,\sigma^*)$-compatible orientation $\overrightarrow{O^{cp}}$ in the cycle-cocycle reversal class of $\overrightarrow{O}$, and $\overrightarrow{O}$ can be  obtained by reversing disjoint directed cycles and cocycles in $\overrightarrow{O^{cp}}$ (see Corollary~\ref{cor0}). Hence the ($\sigma,\sigma^*$)-compatible orientations are representatives for the cycle-cocycle reversal classes. 
Figure~\ref{F-2} gives an acyclic cycle signature $\sigma$ and an acyclic cocycle signature $\sigma^*$ of the triangle graph, which will be fixed for all the examples in this paper. In Figure~\ref{F-1}, the top three orientations are ($\sigma,\sigma^*$)-compatible. For more examples of the two signatures, see \cite[Example 1.1.1 and Example 1.1.3]{BBY}.

Let $T\in\mathcal{T}(G)$ and $e\in E$. If $e\notin T$, then we call the unique cycle in $T\cup \{e\}$ the \emph{fundamental cycle} of $e$ with respect to $T$, denoted by $C(T,e)$; if $e\in T$, then we call the unique cocycle in $(E\backslash T)\cup \{e\}$ the \emph{fundamental cocycle} of $e$ with respect to $T$, denoted by $C^*(T,e)$.

The BBY bijection is as follows. 

\begin{proposition}\label{tree}(\cite{BBY},Theorem 1.3.1(1))
Fix acyclic signatures $\sigma$ and $\sigma^*$ of $G$.  
Then the map $$\text{BBY}_{\sigma,\sigma^*}:\mathcal{T}(G)\longrightarrow\{(\sigma,\sigma^*)\text{-compatible orientations}\}$$ is a bijection, where $\text{BBY}_{\sigma,\sigma^*}$ sends $T$ to the orientation of $G$ in which we orient each $e\notin T$ according to its orientation in $\sigma(C(T,e))$ and each $e\in T$ according to its orientation in $\sigma^*(C^*(T,e))$.
\end{proposition}

Now we construct our extension of the BBY bijection. It is a map from orientations to subgraphs, which actually extends the inverse map $\text{BBY}^{-1}_{\sigma,\sigma^*}$. For an orientation $\overrightarrow{O}$, 
recall that there exists a unique ($\sigma,\sigma^*$)-compatible orientation $\overrightarrow{O^{cp}}$ such that $\overrightarrow{O}$ is obtained from $\overrightarrow{O^{cp}}$ by reversing certain disjoint directed cycles $\{\overrightarrow{C_i}\}_{i\in I}$ and cocycles $\{\overrightarrow{C_j^*}\}_{j\in J}$. By applying the BBY bijection to $\overrightarrow{O^{cp}}$, we get a spanning tree $T=\text{BBY}^{-1}_{\sigma,\sigma^*}(\overrightarrow{O^{cp}})$. Then we add the edges in the reversed cycles to $T$ and delete the edges in the reversed cocycles from $T$, and hence we get a spanning subgraph $S=(T\cup \biguplus\limits_{i\in I}C_i)\backslash \biguplus\limits_{j\in J}C_j^*$. We denote $\overrightarrow{O}\mapsto S$ by $\varphi_{\sigma,\sigma^*}$. Note that if $\overrightarrow{O}$ is ($\sigma,\sigma^*$)-compatible, then no directed cycles or cocycles will be reversed and hence $S=T=\text{BBY}^{-1}_{\sigma,\sigma^*}(\overrightarrow{O})$. Therefore $\varphi_{\sigma,\sigma^*}$ extends $\text{BBY}^{-1}_{\sigma,\sigma^*}$. 

Now we are ready to present the main results of this paper.

\begin{theorem}\label{th}

Fix acyclic signatures $\sigma$ and $\sigma^*$ of $G$. 

(1) The map 

\begin{eqnarray*}
\varphi_{\sigma,\sigma^*}:\{\text{discrete orientations}\} & \longrightarrow & \{\text{spanning subgraphs}\} \\
\overrightarrow{O} & \mapsto & (\text{BBY}_{\sigma,\sigma^*}^{-1}(\overrightarrow{O^{cp}})\cup \biguplus_{i\in I}C_i)\backslash \biguplus_{j\in J}C_j^*
\end{eqnarray*}
is a bijection, where $\overrightarrow{O}$ is an orientation obtained by reversing disjoint directed cycles $\{\overrightarrow{C_i}\}_{i\in I}$ and directed cocycles $\{\overrightarrow{C_j^*}\}_{j\in J}$ in a ($\sigma,\sigma^*$)-compatible orientation $\overrightarrow{O^{cp}}$.

(2) The map $\varphi_{\sigma,\sigma^*}$ specializes to the bijection
\begin{eqnarray*}
\varphi_{\sigma,\sigma^*}: \{\sigma\text{-compatible orientations}\} & \longrightarrow & \{\text{spanning forests}\} \\
\overrightarrow{O} & \mapsto & \text{BBY}_{\sigma,\sigma^*}^{-1}(\overrightarrow{O^{cp}})\backslash \biguplus_{j\in J}C_j^*.
\end{eqnarray*}

(3) The map $\varphi_{\sigma,\sigma^*}$ specializes to the bijection
\begin{eqnarray*}
\varphi_{\sigma,\sigma^*}:\{\sigma^*\text{-compatible orientations}\} & \longrightarrow & \{\text{connected spanning subgraphs}\} \\
\overrightarrow{O} & \mapsto & \text{BBY}_{\sigma,\sigma^*}^{-1}(\overrightarrow{O^{cp}})\cup \biguplus_{i\in I}C_i.
\end{eqnarray*}
\end{theorem}

Note that the $\sigma$-compatible (resp. $\sigma^*$-compatible) orientations form representatives for cycle reversal equivalence classes (resp. cocycle reversal equivalence classes) (see Proposition~\ref{sigmacompatible1}). Hence our construction establishes bijections for all the objects in Table~\ref{table1}.

\begin{example}\label{ex1}
This example illustrates Theorem~\ref{th}. Let $G$ be the triangle graph and two signatures be as in Figure~\ref{F-2}. Figure~\ref{F5} shows the extended BBY bijection $\varphi_{\sigma,\sigma^*}$. It is convenient to draw the orientation $\overrightarrow{O}$ and the subgraph $S$ together provided $\varphi_{\sigma,\sigma^*}(\overrightarrow{O})=S$. The eight orientations are divided into the three cycle-cocycle reversal classes, where the three $(\sigma,\sigma^*)$-compatible orientations correspond to the three spanning trees. One can check that the seven forests correspond to the seven $\sigma$-compatible orientations (see also Figure~\ref{F4}) and the four connected spanning subgraphs correspond to the four $\sigma^*$-compatible orientations (see also Figure~\ref{F7}).

\begin{figure}[h]
            \centering
            \includegraphics[scale=1]{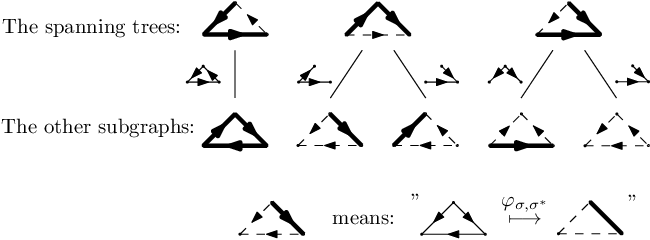}
            \caption{This is the bijection $\varphi_{\sigma,\sigma^*}$ obtained by applying Theorem~\ref{th} to the triangle graph $G$ and the two signatures in Figure~\ref{F-2}. If an orientation is mapped to a subgraph, then we draw them together. An edge is not in the subgraph if and only if it is dashed.}
            \label{F5}
\end{figure}

\end{example}

We mention two other known bijections between spanning subgraphs and orientations with nice specializations. 

Olivier Bernardi defined a family of bijections in \cite{Bernardi}. The bijections rely upon a \emph{combinatorial embedding} of the graph (also known as a \emph{ribbon structure}) and a root \emph{half-edge}. By a result of Chi Ho Yuen \cite[Theorem 20]{Yuen}, for planar graphs, Bernardi's bijection restricted to spanning trees coincides with the BBY bijection for some pair of signatures; see also \cite[Example 1.1.3]{BBY}. However, our extended BBY bijection is different from Bernardi's in general. 

Emeric Gioan and Michel Las Vergnas defined the \emph{active bijections} in \cite{GM}. The active bijections rely upon a \emph{total order} on $E$ and a given reference orientation on $E$. The active bijections provide a canonical bijection between spanning trees and 
\emph{activity classes} of orientations which does not rely upon the reference orientation. The active bijections extend to general oriented matroids \cite{GM2}. The series of papers studying the active bijections started from \cite{G}.

Both of their work are closely related to \emph{activities} and interpret the evaluations of the Tutte polynomial $T_G(x,y)$ for $x=0,1,2$ and $y=0,1,2$. Our work does not involve activities or deal with the case $x=0$ or $y=0$. 

Our bijections rely upon a pair of acyclic signatures which serve as reference orientations on cycles and cocycles. A total order on $E$ together with a given reference orientation on $E$ can induce a pair of acyclic signatures such that every cycle or cocycle in the signatures is oriented according to the reference orientation of its minimum edge \cite[Example 1.1.1]{BBY}. However, the BBY bijection coming from such a pair of signatures is different from the active bijection in general.

\subsection{Introduction to the geometric ideas behind the main theorem with an example}

Now we use an example to illustrate the geometric interpretation of the BBY bijection, which is due to \cite{BBY}, and explain the geometric ideas behind the bijection $\varphi_{\sigma,\sigma^*}$ in Theorem~\ref{th}(2), which will be treated rigorously in Section~\ref{geometric}. The full geometric interpretation of Theorem~\ref{th} will be discussed in Section~\ref{tiling}.

\begin{figure}[h]
            \centering
            \includegraphics[scale=0.9]{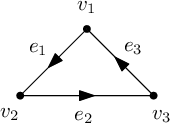}
            \caption{The reference orientation of the triangle graph $G$.}
            \label{F1}
\end{figure}

Let $G$ be the triangle graph with the reference orientation as shown in Figure~\ref{F1}. We identify the continuous orientations of $G$ with the cube $[0,1]^E$, whose vertices are the (discrete) orientations.

The incidence matrix of $G$ is $\begin{pmatrix}
1 & -1 & 0\\
0 & 1 & -1\\
-1 & 0 & 1
\end{pmatrix}$. We delete the last row to get a full-rank matrix $D=\begin{pmatrix}
1 & -1 & 0\\
0 & 1 & -1
\end{pmatrix}$.

\begin{figure}[h]
            \centering
            \includegraphics[scale=0.9]{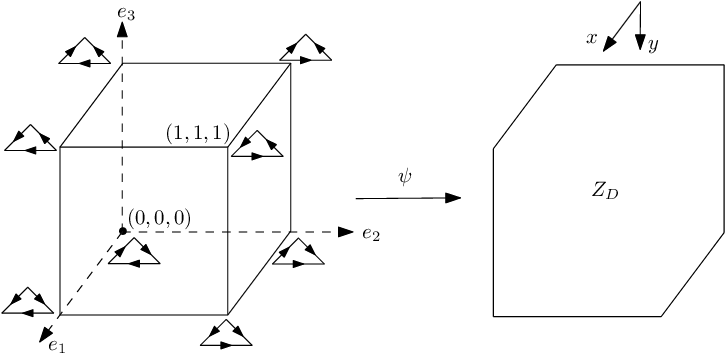}
            \caption{The map $\psi$ and the zonotope $Z_D$. Here we draw $Z_D$ in the indicated coordinate system so that it looks like its section $\im(\mu)$; see Section~\ref{geometric} for the definition of $\im(\mu)$. }
            \label{F12}
\end{figure}

We restrict the linear map $D$ to the continuous orientations $[0,1]^E$. The restricted map is called $\psi$. The map $\psi$ projects the cube $[0,1]^E$ to the \emph{zonotope}\footnote{In general, a zonotope is a \emph{Minkowski sum} of closed line segments and the Minkowski sum is defined by $A+B=\{a+b:a\in A, b\in B\}$, where $A$ and $B$ are two subsets of $\mathbb{R}^n$.} $Z_D=\{\sum\limits_{i=1}^{|E|}c_iv_i:0\leq c_i \leq 1\}$, where $v_i$'s are the columns of $D$; see Figure~\ref{F12}. 

\begin{figure}[h]
            \centering
            \includegraphics[scale=1]{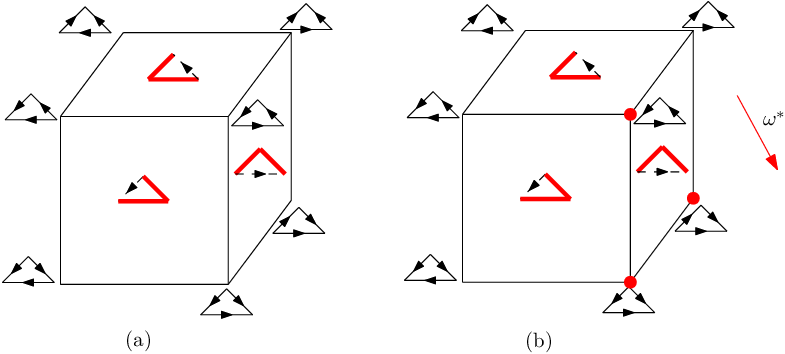}
            \caption{(a): The decomposition $Z_D=\bigcup_{T\in\mathcal{T}(G)}Z_\sigma(T)$. The three spanning trees indicated by solid lines are put in the middle of the corresponding parallelograms. By bi-orienting the edges in the solid lines, one gets the continuous orientations inside the parallelograms. (b): The BBY bijection which sends a spanning tree to the right-bottom vertex of the corresponding parallelogram. }
            \label{F2}
\end{figure}

Fix an acyclic cycle signature $\sigma$. The zonotope $Z_D$ has a polyhedral subdivision
$Z_D=\bigcup_{T\in\mathcal{T}(G)}Z_\sigma(T)$, where $Z_\sigma(T)$ is the image (under $\psi$) of all the continuous orientations where every edge $e\notin T$ is oriented according to the $\sigma$-oriented directed fundamental cycle $\sigma(C(T,e))$; see Figure~\ref{F2}(a). The vertices in the polyhedral subdivision correspond to the (discrete) $\sigma$-compatible orientations and hence by abusing language we may identify each of the vertices with the corresponding orientation. Every tile $Z_\sigma(T)$ is a parallelogram of dimension $|T|$. Each edge of $Z_\sigma(T)$  connects two $\sigma$-compatible orientations, which differ by one arc. The underlying edge of this arc is in the tree $T$. Moreover, the parallel edges of $Z_\sigma(T)$ correspond to the same edge in $T$ and the $|T|$ sets of parallel edges of $Z_\sigma(T)$ correspond to the $|T|$ edges in $T$.

Now we fix a generic vector $w^*$ in the space spanned by $Z_D$ and consider shifting the three parallelograms $Z_\sigma(T)$ along $w^*$. For sufficiently small positive $\epsilon$, the image of $Z_\sigma(T)$ under the shifting map $v\mapsto v+\epsilon w^*$ covers a unique vertex of $Z_\sigma(T)$, which turns out to be the $(\sigma,\sigma^*)$-compatible orientation $\text{BBY}_{\sigma,\sigma^*}(T)$, where $\sigma^*$ is an acyclic cocycle signature determined by $w^*$ and $\text{BBY}_{\sigma,\sigma^*}$ is defined in Proposition~\ref{tree}; see Figure~\ref{F2}(b) for the example and see Figure~\ref{F0}(a) for the general case. So the bijection $\text{BBY}_{\sigma,\sigma^*}$ is the map sending each spanning tree $T$ of $G$ to the $(\sigma,\sigma^*)$-compatible orientation covered by $Z_\sigma(T)+\epsilon w^*$. For proofs of all these statements, see \cite{BBY}.

\begin{figure}[h]
            \centering
            \includegraphics[scale=0.9]{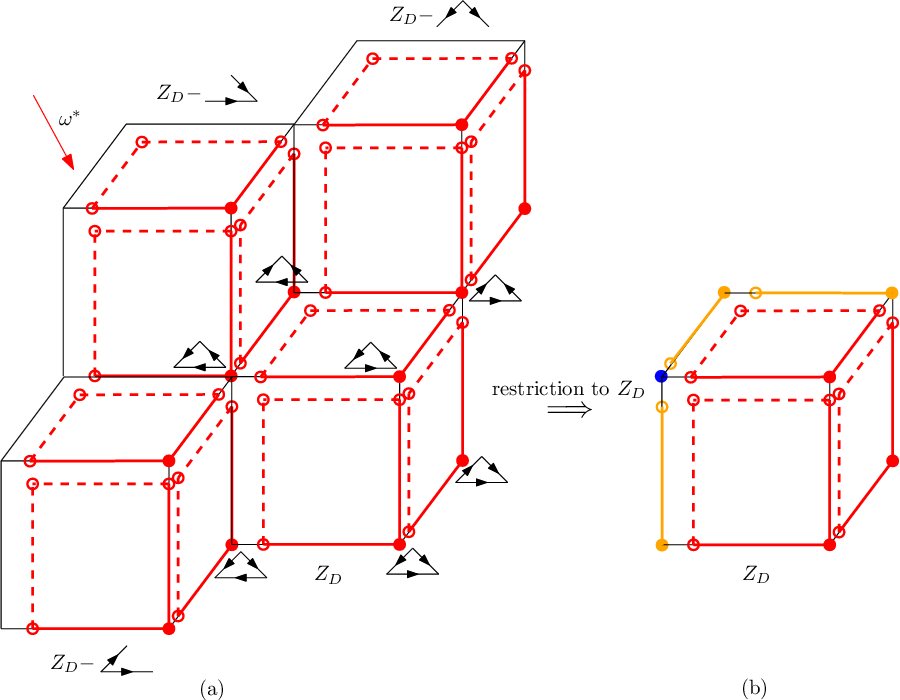}
            \caption{(a) Copies of $Z_D$ tile $\mathbb{R}^2$, although we only draw four copies. Three copies among them are obtained by subtracting three directed cocycles respectively from $Z_D$. The direction $w^*$ induces a half-open decomposition of $\mathbb{R}^2$ indicated by parallelograms in both solid and dashed lines. (b) The tiling of $Z_D$ obtained from restricting the decomposition in (a) to $Z_D$.}
            \label{F3}
\end{figure}

\begin{figure}[h]
            \centering
            \includegraphics[scale=1]{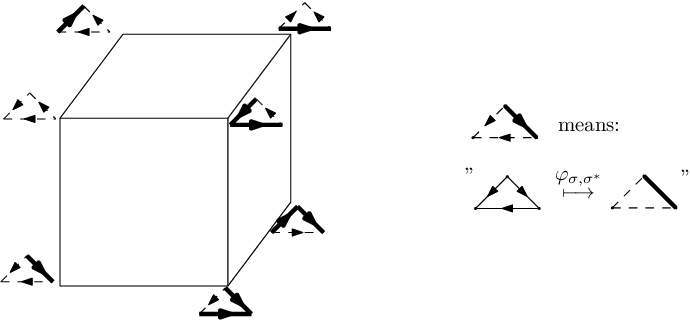}
            \caption{The bijection $\varphi_{\sigma,\sigma^*}$ between spanning forests and $\sigma$-compatible orientations derived from Figure~\ref{F3}(b). For conciseness, if a forest is mapped to an orientation, we combine them into one configuration and put it beside the vertex corresponding to the orientation.}
            \label{F4}
\end{figure}

In order to extend the BBY bijection, we make use of two kinds of tiling related to the zonotope $Z_D$. First we use copies of $Z_D$ to tile its ambient space $\im_\mathbb{R}(D)$, which is $\mathbb{R}^2$ in this example; see Figure~\ref{F3}(a). It can be proved that two adjacent copies of $Z_D$ differ by a sum of disjoint directed cocycles\footnote{Strictly speaking, they differ by $D\cdot \overrightarrow{C^*}$, where $\overrightarrow{C^*}$ is a sum of disjoint directed cocycles (viewed as $\{0,\pm 1\}$-vectors). }.  
The direction $w^*$ induces a decomposition of $\mathbb{R}^2$ into \emph{half-open cells}\footnote{A half-open cell is a Minkowski sum of linear independent half-open line segments.}, each of which contains a unique lattice point. Note that the three half-open cells in the original copy $Z_D$ correspond to the three spanning trees and the lattice points they contain are $(\sigma,\sigma^*)$-compatible orientations. Now we restrict the half-open decomposition of $\mathbb{R}^2$ to $Z_D$ and hence get the second kind of tiling, which is a decomposition of the zonotope into half-open cells; see Figure~\ref{F3}(b). In the tiling, each half-open cell contains a unique lattice point, which corresponds to a $\sigma$-compatible orientation and the vectors that generate the half-open cell correspond to a forest. Hence we get a map from the set of $\sigma$-compatible orientations to the set of forests; see Figure~\ref{F4}. It is easy to see that this map is $\varphi_{\sigma,\sigma^*}$ defined as in Theorem~\ref{th}(2). This geometric approach will be treated rigorously in Section~\ref{geometric}.

\begin{Rem}
In this paper, the two kinds of tiling are not original. For the space tiling zonotopes, see  \cite[Theorem 2]{DG}; for the decomposition of the zonotope into half-open cells, see \cite[Lemma 2.1]{S}. Also see \cite[Section 5.6]{Yuen2} for another application of the space tiling zonotopes in the context of the BBY bijections. 

\end{Rem}

By making use of the duality between cycles and cocycles, we get the map $\varphi_{\sigma,\sigma^*}$ in Theorem~\ref{th}(3). Inspired by the two maps in Theorem~\ref{th}(2) and (3), we were led to conjecture Theorem~\ref{th}(1), but we do not know a geometric proof of it. Instead, we give a purely combinatorial proof. This includes a combinatorial proof of Proposition~\ref{tree}, which was proved in a geometric way in \cite{BBY}. We actually prove a stronger version of Proposition~\ref{tree}, which is used in the proof of Theorem~\ref{th}(1).

\subsection{The structure of the paper}

Now let us outline the structure of the paper. 
In Section~\ref{combinatorial}, we introduce preliminaries and give the combinatorial proofs of Proposition~\ref{tree} and Theorem~\ref{th}. In Section~\ref{geometric}, we briefly review the geometric proof of Proposition~\ref{tree} in \cite{BBY} and give the geometric proof of Theorem~\ref{th}(2).  In Section~\ref{tiling}, we discuss a geometric interpretation of Theorem~\ref{th}, which is a half-open decomposition of the cube $[0,1]^E$. The proof does not make use of zonotopes or translations. We also show that the restriction of the decomposition to the zonotope $Z_D$ (via $\psi$) is exactly the one used to derive Theorem~\ref{th}(2) (e.g. Figure~\ref{F3}(b)). Note that the work of \cite{BBY} is done in the setting of \emph{regular matroids}. While this paper is written mainly in the setting of graphs, we explain in Section~\ref{regular} that the main results hold for regular matroids.

\section{Combinatorial Proof of the Main Results}\label{combinatorial}

In this section, we will prove Proposition~\ref{tree} and Theorem~\ref{th} using combinatorial methods. We will actually prove stronger results, Proposition~\ref{proptree} and Proposition~\ref{propsubset}. We need Proposition~\ref{proptree} to prove Theorem~\ref{th}.  Proposition~\ref{propsubset} will play a key role in Section~\ref{tiling}.

In \cite{BBY}, the proof of Proposition~\ref{tree} is geometric, although the statement is completely combinatorial. Our proof of Proposition~\ref{proptree} is combinatorial, so we obtain a new and combinatorial proof of Proposition~\ref{tree}.



\subsection{Preliminaries}
Let $G$ be a connected finite graph with nonempty edge set $E$, where \emph{loops} and \emph{multiple edges} are allowed.

For each edge $e\in E$, we may assign a direction to it by choosing one of its two endpoints to be the \emph{head} and the other one to be the \emph{tail} and hence get an \emph{arc} (directed from the tail to the head). Note that a loop has two possible directions. An \emph{orientation} of the graph $G$ is an assignment of a direction to each edge, typically denoted by $\overrightarrow{O}$. 
An \emph{partial orientation} of the graph $G$ is an assignment of a direction to each edge in a subset of $E$, typically denoted by $\overrightarrow{P}$.

A subset $C$ of $E$ is called a \emph{cycle} if there exists distinct vertices $v_1,v_2, \cdots, v_n$ such that $C=\{\text{edge }v_iv_{i+1}:i=1,2,\cdots,n\}$, 
where $v_{n+1}:=v_1$. Note that a cycle may be a loop. If we direct every edge in $C$ from $v_i$ to $v_{i+1}$ or direct every edge in $C$ from $v_{i+1}$ to $v_{i}$, then we get a \emph{directed cycle}, which is typically denoted by $\overrightarrow{C}$. Given a subset $W$ of vertices, the set of edges with one endpoint in $W$ and the other one not in $W$ is called a \emph{cut}. A \emph{cocycle} $C^*$ is a cut which is minimal for inclusion (equivalently it is a cut whose deletion increases the number of connected components by one). If we direct every edge in $C^*$ from $W$ to its complement or directed every edge in the other way, then we get a \emph{directed cocycle}, which is typically denoted by $\overrightarrow{C^*}$. 

When an arc $\overrightarrow{e}$, a directed cycle $\overrightarrow{C}$, a directed cocycle $\overrightarrow{C^*}$, or a partial orientation $\overrightarrow{P}$ is specified, the corresponding underlying edge(s) will be denoted by $e$, $C$, $C^*$, or $P$, respectively. Viewing $\overrightarrow{O}$, $\overrightarrow{C}$, $\overrightarrow{C^*}$, and $\overrightarrow{P}$ as sets of arcs, it makes sense to write $\overrightarrow{e}\in\overrightarrow{O}$, etc. When $\overrightarrow{P}$ is a subset of $\overrightarrow{O}$, we say $\overrightarrow{P}$ is \emph{in} the orientation $\overrightarrow{O}$. In particular, $\overrightarrow{P}$ can be a directed cycle or a directed cocycle.

Given a spanning tree $T$ and an arc $\overrightarrow{e}$, we denote by $C(T,\overrightarrow{e})$ the fundamental cycle directed according to $\overrightarrow{e}$ when $e\in T$, and denote by 
$C^*(T,\overrightarrow{e})$ the fundamental cocycle directed according to $\overrightarrow{e}$ when $e\notin T$.

It is a classical fact that every arc in an orientation belongs to either a directed cycle or a directed cocycle (but not both). Here we need the following stronger version.

\begin{lemma}\label{3-painting}
Let $E_c$ and $E_d$ be two disjoint subsets of $E$, $\overrightarrow{P}$ be a partial orientation with support $E\backslash (E_c\cup E_d)$, and $\overrightarrow{e}$ be an arc in $\overrightarrow{P}$. Then there exists a directed cycle $\overrightarrow{C}$ containing $\overrightarrow{e}$ such that $\overrightarrow{C}$ agrees with $\overrightarrow{P}$ on any edge they share and $C\cap E_d=\emptyset$, or there exists a directed cocycle $\overrightarrow{C^*}$ containing $\overrightarrow{e}$ such that $\overrightarrow{C^*}$ agrees with $\overrightarrow{P}$ on any edge they share and $C^*\cap E_c=\emptyset$. 
\end{lemma}

\begin{proof}
Assume the arc $\overrightarrow{e}$ is directed from the vertex $u$ to $w$. Let $W$ be set of vertices reachable from $w$, where one may use the arcs in $\overrightarrow{P}$ and both directions of the edges in $E_c$. If $u\in W$, then there exists a directed cycle $\overrightarrow{C}$ as desired. Otherwise $\overrightarrow{e}$ belongs to the directed cut $\overrightarrow{C^*_0}$ oriented from the complement of $W$ to $W$. Then the directed cocycle (minimal directed cut) $\overrightarrow{C^*}$ containing $\overrightarrow{e}$ in $\overrightarrow{C^*_0}$ is as desired. 
\end{proof}

\begin{Rem}
(1) We do not need the fact that the two cases in Lemma~\ref{3-painting} are exclusive. For the proof of the exclusiveness, see either of the two references below. 

(2) Lemma~\ref{3-painting} (including the proof) is essentially the same as Proposition 2.5 in \cite{BH}, which is a paper studying \emph{fourientations} of a graph.

(3) Lemma~\ref{3-painting} can be generalized to \emph{regular matroids} or even \emph{oriented matroids}, and it is called the \emph{3-painting axiom}; see \cite[Theorem 3.4.4]{BVSWZ}. However, our proof cannot be generalized accordingly, because it makes use of vertices. 
\end{Rem}

By fixing a reference orientation, we identify the set of orientations of $G$ with the set $\{0,1\}^E$, where the point $(1,1,\cdots,1)$ corresponds to the reference orientation. Let $D$ be the modified incidence matrix of the reference orientation with the last row removed. To be precise, if we denote the vertices of $G$ by $v_1, v_2, \cdots, v_{r+1}$ and the arcs in the reference orientation by $\overrightarrow{e_1}, \overrightarrow{e_2}, \cdots, \overrightarrow{e_n}$, then $D$ is the $r \times n$ matrix $(d_{ij})$ whose entries are 

\begin{equation*}
  d_{ij}=
  \begin{cases}
    1, &\text{if vertex $v_i$ is the head of non-loop arc $\overrightarrow{e_j}$;} \\
    -1, &\text{if vertex $v_i$ is the tail of non-loop arc $\overrightarrow{e_j}$;}\\
    0, & \text{otherwise}.
  \end{cases}
\end{equation*}

See Section~\ref{intro} for an example of $D$. 

One minor issue here is that when the graph $G$ only has one vertex ($r=0$), the matrix $D$ is not defined. In this (trivial) case, it is easy to check that Theorem~\ref{th} holds. Note that some geometric results do not hold, whether we use the incidence matrix or the matrix $D$. For example, Proposition~\ref{bbymain} does not hold because $Z_\sigma(T)$ is a single point if we use the incidence matrix and is not defined if we use the matrix $D$. From now on, we assume $G$ contains at least two vertices. 

\begin{Rem}
The advantage of removing one row is to make $D$ of full rank. Hence in Section~\ref{geometric}, the zonotope $Z_D$ spans the whole space rather than a hyperplane. It does not matter which row is removed. Also note that when we use a matrix to represent a regular matroid in Section~\ref{regular}, the matrix is chosen to be of full rank.
\end{Rem}

Now we give some general facts concerning the directed cycles and cocycles in terms of linear algebra. We refer readers to \cite{Biggs} for the details. The vector space $\mathbb{R}^E$ has an orthogonal decomposition $\ker_\mathbb{R}(D)\oplus\im_\mathbb{R}(D^T)$ with respect to the standard inner product $\langle\cdot,\cdot\rangle$. We call $\ker_\mathbb{R}(D)$ the \emph{cycle space} and  $\im_\mathbb{R}(D^T)$ the \emph{cocycle space}. Any directed cycle (resp. directed cocycle), written as a $\{0,\pm1\}$-vector according the reference orientation, is an element in the cycle space (resp. cocycle space). For lattice points, we denote the two abelian groups $\ker_\mathbb{Z}(D)=\ker_\mathbb{R}(D)\cap \mathbb{Z}^E$ and $\im_\mathbb{Z}(D^T)=\im_\mathbb{R}(D^T)\cap \mathbb{Z}^E$. For any spanning tree $T$, the directed fundamental cycles (resp. directed fundamental cocycles) form a basis of $\ker_\mathbb{R}(D)$ (resp. $\im_\mathbb{R}(D^T)$), and an integral basis of $\ker_\mathbb{Z}(D)$ (resp. $\im_\mathbb{Z}(D^T)$). 

We recall the \emph{cycle reversal (equivalence) classes}, \emph{cocycle reversal (equivalence) classes}, and \emph{cycle-cocycle reversal (equivalence) classes} of orientations of $G$ introduced in \cite{G1}. If $\overrightarrow{C}$ is a directed cycle in $\overrightarrow{O}$, then a \emph{cycle reversal} replaces each arc in $\overrightarrow{C}$ with the opposite arc, and hence gives a new orientation, which is $\overrightarrow{O}-\overrightarrow{C}$ written as a vector, where $\overrightarrow{O}$ is a \{0,1\}-vector and $\overrightarrow{C}$ is a $\{0,\pm1\}$-vector. The equivalence relation generated by cycle reversals defines the cycle reversal classes of orientations of $G$. Similarly, we define the cocycle reversal classes. The equivalence relation generated by cycle and cocycle reversals defines the cycle-cocycle reversal classes. 

When $\overrightarrow{P}$ is \emph{in} the orientation $\overrightarrow{O}$, we denote by $-_{\overrightarrow{P}}\overrightarrow{O}$ or $-_{P}\overrightarrow{O}$ the orientation obtained by reversing $\overrightarrow{P}$ in $\overrightarrow{O}$. For example, if $\overrightarrow{C}$ is a directed cycle in $\overrightarrow{O}$, then $\overrightarrow{O}-\overrightarrow{C}=-_{\overrightarrow{C}}\overrightarrow{O}$.

We need the following lemmas. 

\begin{lemma}\cite[Lemma 6.7]{Z}\label{conformal0}
(1) Let $\overrightarrow{u}\in\mathbb{R}^E$ be a vector in $\ker_\mathbb{R}(D)$. Then $\overrightarrow{u}$ can be written as a sum of directed cycles with positive coefficients $\sum k_i\overrightarrow{C_i}$ where for each edge $e$ of each $C_i$, the sign of $e$ in $\overrightarrow{C_i}$ agrees with the sign of $e$ in $\overrightarrow{u}$. 

(2) Let $\overrightarrow{u}\in\mathbb{R}^E$ be a vector in $\im_\mathbb{R}(D^T)$. Then $\overrightarrow{u}$ can be written as a sum of directed cocycles with positive coefficients $\sum k_i\overrightarrow{C_i^*}$ where for each edge $e$ of each $C_i^*$, the sign of $e$ in $\overrightarrow{C_i^*}$ agrees with the sign of $e$ in $\overrightarrow{u}$. 
\end{lemma}

Sometimes the sum of directed cycles and the sum of directed cocycles are also denoted by $\overrightarrow{C}$ and $\overrightarrow{C^*}$, respectively. 
\begin{lemma}\label{conformal}

(1) If $\overrightarrow{C}\in\ker_\mathbb{Z}(D)$ is a $\{0,\pm1\}$-vector, then $\overrightarrow{C}$ is a sum of disjoint directed cycles. 

(2) If $\overrightarrow{C^*}\in\im_\mathbb{Z}(D^T)$ is a $\{0,\pm1\}$-vector, then $\overrightarrow{C^*}$ is a sum of disjoint directed cocycles. 

\end{lemma}

\begin{proof}
(1) This is from \cite[Lemma 4.1.1]{BBY}.

(2) The proof of \cite[Lemma 4.1.1]{BBY} is derived from Lemma~\ref{conformal0}(or \cite[Lemma 6.7]{Z}), so the dual argument works.  
\end{proof}

The next lemma gives two simple descriptions of two orientations being in the same cycle reversal class. 

\begin{lemma}\label{cyclelemma}
Let $\overrightarrow{O_1}$ and $\overrightarrow{O_2}$ be two orientations of $G$. Then the following are equivalent:

(a) $\overrightarrow{O_1}$ and $\overrightarrow{O_2}$ are in the same cycle reversal class;

(b) $\overrightarrow{O_1}-\overrightarrow{O_2}\in\ker_\mathbb{Z}(D)$;

(c) $\overrightarrow{O_1}$ can be obtained from $\overrightarrow{O_2}$ by  reversing disjoint directed cycles. 

\end{lemma}
\begin{proof}
It is trivial that (a) implies (b) and (c) implies (a). By Lemma~\ref{conformal}, (b) implies (c). 
\end{proof}

Similar results hold for cocycles.

Recall the definitions of the acyclic cycle signature $\sigma$, the acyclic cocycle signature $\sigma^*$, the $\sigma$-compatible orientations, the $\sigma^*$-compatible orientations, and the $(\sigma,\sigma^*)$-compatible orientations in Section~\ref{intro}. The next result shows that $\sigma$-compatible orientations are representatives of cycle reversal classes.  

\begin{proposition}\cite[Prop. 4.1.4]{BBY}\label{sigmacompatible1}
Let $\sigma$ be an acyclic cycle signature. Then each cycle reversal class contains a unique $\sigma$-compatible orientation.  
\end{proposition}

The similar result to Proposition~\ref{sigmacompatible1} holds for cocycles.

\begin{corollary}\label{cor0}
Let $\sigma$ be an acyclic cycle signature,  $\sigma^*$ be an acyclic cocycle signature, and $\overrightarrow{O}$ be an orientation.

(1) There exists a unique $\sigma$-compatible orientation $\overrightarrow{O^\sigma}$ in the cycle reversal class of $\overrightarrow{O}$, and $\overrightarrow{O}$ can be obtained by reversing disjoint directed cycles in $\overrightarrow{O^\sigma}$.

(2) There exists a unique $\sigma^*$-compatible orientation $\overrightarrow{O^{\sigma^*}}$ in the cocycle reversal class of $\overrightarrow{O}$, and $\overrightarrow{O}$ can be obtained by reversing disjoint directed cocycles in $\overrightarrow{O^{\sigma^*}}$.

(3) There exists a unique $(\sigma,\sigma^*)$-compatible orientation $\overrightarrow{O^{cp}}$ in the cycle-cocycle reversal class of $\overrightarrow{O}$, and $\overrightarrow{O}$ can be obtained by reversing disjoint directed cycles and cocycles in $\overrightarrow{O^{cp}}$. Moreover, if $\overrightarrow{O}$ is already $\sigma$-compatible, then only directed cocycles are reversed; if $\overrightarrow{O}$ is already $\sigma^*$-compatible, then only directed cycles are reversed. 

\end{corollary}
\begin{proof}
By Proposition~\ref{sigmacompatible1} and Lemma~\ref{cyclelemma}, we have (1). (2) is dual to (1). Because every arc in an orientation belongs to a directed cycle or a directed cocycle but not both, (3) holds.
\end{proof}

\subsection{Combinatorial proof of Proposition~\ref{tree}}

We fix a graph $G$, a reference orientation, an acyclic cycle signature $\sigma$, and an acyclic cocycle signature $\sigma^*$.

Let $e$ be an edge. For an orientation $\overrightarrow{O}$, we denote by $\overrightarrow{O}(e)$ the arc $\overrightarrow{e}$ (whose underlying edge is $e$) in $\overrightarrow{O}$. If $e$ is in the support of a directed cycle $\overrightarrow{C}$ (resp. directed cocycle $\overrightarrow{C^*}$), we denote by $\overrightarrow{C}(e)$ (resp. $\overrightarrow{C^*}(e)$) the arc $\overrightarrow{e}$ in the directed cycle $\overrightarrow{C}$ (resp. directed cocycle $\overrightarrow{C^*}$). More generally, by viewing $\overrightarrow{v}\in\mathbb{Z}^E$ as a multiset of arcs, we denote by $\overrightarrow{v}(e)$ the arc $\overrightarrow{e}$ (with multiplicity one) in $\overrightarrow{v}$ if $e$ in the support of $\overrightarrow{v}$. In this case, we also write  $\overrightarrow{e}\in\overrightarrow{v}$. 

The following lemma plays an important role in our proofs.
\begin{lemma}\label{fundamental}
Fix a spanning tree $T$. 

(1) For any $\overrightarrow{C}\in\ker_\mathbb{Z}(D)$, $$\overrightarrow{C}=\sum_{e\notin T,\overrightarrow{e}\in\overrightarrow{C}}C(T,\overrightarrow{e}),$$

(2) and for any $\overrightarrow{C^*}\in\im_\mathbb{Z}(D^T)$, $$\overrightarrow{C^*}=\sum_{e\in T,\overrightarrow{e}\in\overrightarrow{C^*}}C^*(T,\overrightarrow{e}),$$ where $\overrightarrow{e}$ appears in the summations as many times as its multiplicity. 
\end{lemma}
\begin{proof}
(1) Because the directed fundamental cycles form a basis of $\ker_\mathbb{R}(D)$, we can write $\overrightarrow{C}$ as a linear combination of them.  Note that each edge $e\notin T$ is in exactly one fundamental cycle, so by comparing the coefficients of $e$ in both sides, we get the desired formula. 

(2) The proof is similar. 
\end{proof}

Recall that the map $\text{BBY}_{\sigma,\sigma^*}$ defined in Proposition~\ref{tree} sends a spanning tree $T\in\mathcal{T}(G)$ to the orientation $\overrightarrow{O}$ where $C(T,\overrightarrow{e})\in\sigma$ or $C^*(T,\overrightarrow{e})\in\sigma^*$
for any arc $\overrightarrow{e}\in\overrightarrow{O}$. The next lemma is a key tool in this subsection. It uses the information carried by $\text{BBY}_{\sigma,\sigma^*}(T)$ rather than the two signatures to give a sufficient condition of directed cycles (resp. directed cocycles) being $\sigma$-compatible (resp. $\sigma^*$-compatible). 

\begin{lemma}\label{mainlemma}
Assume $\text{BBY}_{\sigma,\sigma^*}(T)=\overrightarrow{O}$, where $T$ is a spanning tree. 

(1) Let $\overrightarrow{C}\in\ker_\mathbb{Z}(D)$. If for any arc $\overrightarrow{e}\in\overrightarrow{C}$ such that $e\notin T$, we have $\overrightarrow{e}\in\overrightarrow{O}$, then $\overrightarrow{C}$ is a sum of $\sigma$-compatible directed cycles.  In particular, if $\overrightarrow{C}$ is a directed cycle satisfying the condition, then $\overrightarrow{C}$ is $\sigma$-compatible. 

(2) Let $\overrightarrow{C^*}\in\im_\mathbb{Z}(D^T)$. If for any arc $\overrightarrow{e}\in\overrightarrow{C^*}$ such that $e\in T$, we have $\overrightarrow{e}\in\overrightarrow{O}$, then $\overrightarrow{C^*}$ is a sum of $\sigma^*$-compatible directed cocycles. In particular, if $\overrightarrow{C^*}$ is a directed cocycle satisfying the condition, then $\overrightarrow{C^*}$ is $\sigma^*$-compatible.
\end{lemma}

\begin{proof}
(1) By Lemma~\ref{fundamental}, $\overrightarrow{C}=\sum\limits_{e\notin T,\overrightarrow{e}\in\overrightarrow{C}}C(T,\overrightarrow{e})$. By assumption, $\overrightarrow{e}\in\overrightarrow{O}$ and hence the directed cycle $C(T,\overrightarrow{e})$ is $\sigma$-compatible by the definition of $\text{BBY}_{\sigma,\sigma^*}$. 

If $\overrightarrow{C}$ is a directed cycle, then $\overrightarrow{C}$ is $\sigma$-compatible because the signature $\sigma$ is acyclic. 

(2) The proof is similar. 
\end{proof}

We now show that $\text{BBY}_{\sigma,\sigma^*}$ is a map from $\mathcal{T}(G)$ to
$\{(\sigma,\sigma^*)\text{-compatible orientations}\}$.

\begin{lemma}
For any $T\in\mathcal{T}(G)$, the orientation $\text{BBY}_{\sigma,\sigma^*}(T)$ is $(\sigma,\sigma^*)$-compatible.
\end{lemma}
\begin{proof}
By definition, we need to show any directed cycle $\overrightarrow{C}$ in $\text{BBY}_{\sigma,\sigma^*}(T)$ is $\sigma$-compatible 
and any directed cocycle $\overrightarrow{C^*}$ in $\text{BBY}_{\sigma,\sigma^*}(T)$ is $\sigma^*$-compatible. This is a direct consequence of Lemma~\ref{mainlemma}.
\end{proof}

Next we show that $\text{BBY}_{\sigma,\sigma^*}$ is injective. The following result is actually stronger than injectivity, which will be useful later.

\begin{proposition}\label{proptree}
Let $\overrightarrow{O_1}=\text{BBY}_{\sigma,\sigma^*}(T_1)$ and $\overrightarrow{O_2}=\text{BBY}_{\sigma,\sigma^*}(T_2)$, where $T_1$ and $T_2$ are two 
different spanning trees. Then there exists an edge $e$ such that $e\in T_1\bigtriangleup T_2$ and $\overrightarrow{O_1}(e)\neq \overrightarrow{O_2}(e)$. In particular, $\text{BBY}_{\sigma,\sigma^*}$ is injective. 
\end{proposition}

\begin{proof}
Denote $E_c=T_1\cap T_2$ and $E_d=E\backslash(T_1\cup T_2)$. Assume by contradiction that there is no such edge, which means $\overrightarrow{O_1}$ and $\overrightarrow{O_2}$ agree on $E\backslash (E_c\cup E_d)$. 

Let $\overrightarrow{P}$ be the (non-empty) partial orientation obtained by restricting the orientation $-_{T_1}\overrightarrow{O_1}$ to $E\backslash (E_c\cup E_d)$. Note that $\overrightarrow{P}$ is also the restriction of $-_{T_1}\overrightarrow{O_2}$ to $E\backslash (E_c\cup E_d)$.

Apply Lemma~\ref{3-painting} to any arc in $\overrightarrow{P}$. Then one of the following two cases occurs. 

Case 1: There is a directed cycle $\overrightarrow{C}$ such that it agrees with $\overrightarrow{P}$ on any edge they share and $C\cap E_d=\emptyset$. 

Let $\overrightarrow{e}$ be an arc of $\overrightarrow{C}$. Then $e$ is in $T_1\backslash T_2$, $T_2\backslash T_1$ or $E_c$. When $e\in T_1\backslash T_2$, which is equivalent to $e\in C\backslash T_2$, $\overrightarrow{e}=-\overrightarrow{O_2}(e)$; when $e\in T_2\backslash T_1$, which is equivalent to $e\in C\backslash T_1$, $\overrightarrow{e}=\overrightarrow{O_1}(e)$.

Applying Lemma~\ref{mainlemma} to $T_1$ and $\overrightarrow{C}$, we get $\overrightarrow{C}$ is $\sigma$-compatible. 

Applying Lemma~\ref{mainlemma} to $T_2$ and $-\overrightarrow{C}$, we get $-\overrightarrow{C}$ is $\sigma$-compatible. Hence we get a contradiction. 

Case 2: There is a directed cocycle $\overrightarrow{C^*}$ such that it agrees with $\overrightarrow{P}$ on any edge they share and $C^*\cap E_c=\emptyset$. 

This case is dual to Case 1 and hence ends with a contradiction. 

\end{proof}

By Corollary~\ref{cor0}(3), $(\sigma,\sigma^*)$-compatible orientations are representatives of the cycle-cocycle reversal classes. By \cite[Corollary 4.13]{G1}, the number of the cycle-cocycle reversal classes is $T_G(1,1)$ and hence is equal to the number of the spanning trees. Note that the counterpart of this enumerative fact also holds for regular matroids; see \cite[Theorem 3.10]{G2}. So we get the bijectivity.

\begin{corollary}\cite{BBY}
The map $$\text{BBY}_{\sigma,\sigma^*}:\mathcal{T}(G)\longrightarrow\{(\sigma,\sigma^*)\text{-compatible orientations}\}$$ is a bijection. 
\end{corollary}

This concludes the new combinatorial proof of Proposition~\ref{tree}.

\subsection{Combinatorial proof of Theorem~\ref{th}}
We now prove our main result, Theorem~\ref{th}, and a ``local bijectivity'' property. 

We extend $\text{BBY}_{\sigma,\sigma^*}^{-1}$ to $$\varphi_{\sigma,\sigma^*}:\{\text{discrete orientations}\}\longrightarrow \{\text{spanning subgraphs}\}.$$
By Corollary~\ref{cor0}, for any orientation $\overrightarrow{O}$, there exists a unique $(\sigma,\sigma^*)$-compatible orientation $\overrightarrow{O^{cp}}$ such that $\overrightarrow{O}=-_{\overrightarrow{P}}\overrightarrow{O^{cp}}$, where $\overrightarrow{P}$ is the disjoint union of certain directed cycles $\{\overrightarrow{C_i}\}_{i\in I}$ and directed cocycles $\{\overrightarrow{C^*_j}\}_{j\in J}$. Define $$\varphi_{\sigma,\sigma^*}(\overrightarrow{O})=(\text{BBY}_{\sigma,\sigma^*}^{-1}(\overrightarrow{O^{cp}})\cup\biguplus_{i\in I}C_i)\backslash\biguplus_{j\in J}C^*_j.$$

Note that $\varphi_{\sigma,\sigma^*}=\text{BBY}_{\sigma,\sigma^*}^{-1}$ when restricted to the set  $\{(\sigma,\sigma^*)\text{-compatible orientations}\}$. The domain and codomain of  $\varphi_{\sigma,\sigma^*}$ clearly have the same cardinality, so it is enough to show the injectivity. As in Proposition~\ref{proptree}, we prove a stronger property, which will be related to the half-open decomposition in  Section~\ref{tiling}.

\begin{proposition}\label{propsubset}
Let $\overrightarrow{O_1}$ and $\overrightarrow{O_2}$ be two different orientations of $G$. Then there exists an edge $e$ such that $\overrightarrow{O_1}(e)\neq\overrightarrow{O_2}(e)$ and $e\in \varphi_{\sigma,\sigma^*}(\overrightarrow{O_1})\bigtriangleup \varphi_{\sigma,\sigma^*}(\overrightarrow{O_2})$. In particular, $\varphi_{\sigma,\sigma^*}$ is injective. 
\end{proposition}

\begin{proof}

Denote $E_{\rightrightarrows}=\{e:\overrightarrow{O_1}(e)=\overrightarrow{O_2}(e)\}$ and $E_{\rightleftarrows}=\{e:\overrightarrow{O_1}(e)\neq\overrightarrow{O_2}(e)\}$($\neq\emptyset$).

For $k=1,2$, we denote $-_{\overrightarrow{P_k}}\overrightarrow{O_k}=\overrightarrow{O_k^{cp}}$ and $T_k=\varphi_{\sigma,\sigma^*}(\overrightarrow{O_k^{cp}})$, where $\overrightarrow{O_k^{cp}}$ is a $(\sigma,\sigma^*)$-compatible orientation, and $\overrightarrow{P_k}$ is a partial orientation consisting of certain disjoint directed cycles $\{\overrightarrow{C_{k,i}}\}_{i\in I_k}$ and directed cocycles $\{\overrightarrow{C^*_{k,j}}\}_{j\in J_k}$ in  $\overrightarrow{O_k}$.

The proof consists of two parts. 

For the first part, we claim if one of the following cases happens, we get the desired edge. 

Case (a): There exists an edge $e\in E_{\rightleftarrows}\backslash T_1$ such that $e\in (\biguplus_{i\in I_2} C_{2,i})\backslash(\biguplus_{i\in I_1} C_{1,i})$.

Case (b): There exists an edge $e\in E_{\rightleftarrows}\backslash T_2$ such that $e\in (\biguplus_{i\in I_1} C_{1,i})\backslash(\biguplus_{i\in I_2} C_{2,i})$.

Case (a*): There exists an edge $e\in E_{\rightleftarrows}\cap T_1$ such that $e\in (\biguplus_{j\in J_2} C^*_{2,j})\backslash(\biguplus_{j\in J_1} C^*_{1,j})$.

Case (b*): There exists an edge $e\in E_{\rightleftarrows}\cap T_2$ such that $e\in (\biguplus_{j\in J_1} C^*_{1,j})\backslash(\biguplus_{j\in J_2} C^*_{2,j})$.

Case (c): $(\biguplus_{i\in I_2}C_{2,i})\cap(\biguplus_{j\in J_1}C^*_{1,j})\neq\emptyset$.

Case (c*): $(\biguplus_{i\in I_1}C_{1,i})\cap(\biguplus_{j\in J_2}C^*_{2,j})\neq\emptyset$.

If Case (a) holds, then $e\notin\varphi_{\sigma,\sigma^*}(\overrightarrow{O_1})$ and $e\in\varphi_{\sigma,\sigma^*}(\overrightarrow{O_2})$ by the definition of $\varphi_{\sigma,\sigma^*}$. Hence we obtain the desired edge. 

By exchanging the subscripts $1$ and $2$ or taking the dual case, we also obtain the desired edge in Case (b), (a*), and (b*). 

If Case (c) holds, then $C_{2,i}\cap C^*_{1,j}\neq\emptyset$ for some $i,j$. The intersection contains at least two edges and one of the edges, called $e$, is oriented in opposite ways by $\overrightarrow{C_{2,i}}$ and $\overrightarrow{C^*_{1,j}}$. Then $e\in E_{\rightleftarrows}$,  $e\in\varphi_{\sigma,\sigma^*}(\overrightarrow{O_2})$, and $e\notin\varphi_{\sigma,\sigma^*}(\overrightarrow{O_1})$.

By exchanging the subscripts $1$ and $2$, we get Case (c*).

For the second part of the proof, we consider $$\overrightarrow{C}:=\sum_{i\in I_2}\overrightarrow{C_{2,i}}-\sum_{i\in I_1}\overrightarrow{C_{1,i}}\in\ker_\mathbb{Z}(A)$$ and $$\overrightarrow{C^*}:=\sum_{j\in J_2}\overrightarrow{C^*_{2,j}}-\sum_{j\in J_1}\overrightarrow{C^*_{1,j}}\in\im_\mathbb{Z}(A^T).$$ Both of them are $\{0,\pm 1, \pm 2\}$-vectors.

Now we claim that if none of the cases in the first part of the proof holds, then $\overrightarrow{C}=\overrightarrow{C^*}=0$. 

Assume by contradiction that $\overrightarrow{C}\neq 0$.

We consider applying Lemma~\ref{mainlemma} to $T_1$ and $\overrightarrow{C}$. For any arc $\overrightarrow{e}\in\overrightarrow{C}$ such that $e\notin T_1$, we will show $\overrightarrow{e}\in\overrightarrow{O_1^{cp}}$ and hence $\overrightarrow{C}$ is a sum of $\sigma$-compatible directed cycles. 

If the underlying edge of this arc $e\in\biguplus C_{1,i}$, then $\overrightarrow{e}=-(\biguplus \overrightarrow{C_{1,i}})(e)\in\overrightarrow{O_1^{cp}}$. 

If $e\notin\biguplus C_{1,i}$, then $e\in\biguplus C_{2,i}$. Moreover, $e$ must be in $E_{\rightrightarrows}$ because otherwise $e\in E_{\rightleftarrows}$ and hence Case (a) holds. 
Hence $\overrightarrow{e}=\overrightarrow{O_2}(e)=\overrightarrow{O_1}(e)$. Because Case (c) does not hold, $e\in \biguplus C_{2,i}$ implies $e\notin \biguplus C^*_{1,j}$. By $e\notin \biguplus C^*_{1,j}$ and $e\notin \biguplus C_{1,i}$, we get $\overrightarrow{O_1}(e)=\overrightarrow{O_1^{cp}}(e)$. So $\overrightarrow{e}=\overrightarrow{O_1^{cp}}(e)\in\overrightarrow{O_1^{cp}}$. 

Therefore $\overrightarrow{C}$ is a sum of $\sigma$-compatible directed cycles. Similarly, we apply Lemma~\ref{mainlemma} to $T_2$ and $-\overrightarrow{C}=\sum\limits_{i\in I_1}\overrightarrow{C_{1,i}}-\sum\limits_{i\in I_2}\overrightarrow{C_{2,i}}$. Then we get $-\overrightarrow{C}$ is also a sum of $\sigma$-compatible directed cycles, which contradicts that  the cycle signature $\sigma$ is acyclic. 

Therefore, $\overrightarrow{C}=0$.

By a dual argument, $\overrightarrow{C^*}=0$.

It remains to show the proposition is true when $\overrightarrow{C}=\overrightarrow{C^*}=0$.

In this case, the cycles $C_{1,i}$ and $C_{2,i}$ and the cocycles $C^*_{1,j}$ and $C^*_{2,j}$ are all subsets of $E_{\rightrightarrows}$ and hence $\{e:\overrightarrow{O_1^{cp}}(e)\neq\overrightarrow{O_2^{cp}}(e)\}=E_{\rightleftarrows}$. By Proposition~\ref{proptree}, there exists an edge $e\in E_{\rightleftarrows}$ such that $e\in T_1\bigtriangleup T_2$. Since $e\notin(\biguplus_{i\in I_1}C_{1,i})\cup(\biguplus_{j\in J_1}C^*_{1,j})$, $e\in \varphi_{\sigma,\sigma^*}(\overrightarrow{O_1})\bigtriangleup \varphi_{\sigma,\sigma^*}(\overrightarrow{O_2})$ as desired.

\end{proof}

As a consequence, Theorem~\ref{th}(1) holds. 

\begin{Rem}
In Proposition~\ref{propsubset}, if we only want to prove the injectivity of $\varphi_{\sigma,\sigma^*}$ (Theorem~\ref{th}(1)), we will still need Proposition~\ref{proptree}. The bijectivity of $\text{BBY}_{\sigma,\sigma^*}$ is not strong enough to finish the proof.
\end{Rem}

Now we study some specializations of the bijection $\varphi_{\sigma,\sigma^*}$. Theorem~\ref{th}(2) and (3) hold due to the following lemma.

\begin{lemma}\label{l1}
Let $\overrightarrow{O}$ be an orientation and $\overrightarrow{O^{cp}}$ be the $(\sigma,\sigma^*)$-compatible orientation such that one gets $\overrightarrow{O}$ by reversing disjoint directed cycles $\{\overrightarrow{C_i}\}_{i\in I}$ and directed cocycles $\{\overrightarrow{C^*_j}\}_{j\in J}$ in $\overrightarrow{O^{cp}}$. Then

(1) $\varphi_{\sigma,\sigma^*}(\overrightarrow{O})$ is a forest $\Leftrightarrow$ $I=\emptyset$ $\Leftrightarrow$ $\overrightarrow{O}$ is $\sigma$-compatible;

(2) $\varphi_{\sigma,\sigma^*}(\overrightarrow{O})$ is a connected spanning subgraph $\Leftrightarrow$ $J=\emptyset$ $\Leftrightarrow$ $\overrightarrow{O}$ is $\sigma^*$-compatible.
\end{lemma}

\begin{proof}
(1) Recall that $$\varphi_{\sigma,\sigma^*}(\overrightarrow{O})=(\text{BBY}_{\sigma,\sigma^*}^{-1}(\overrightarrow{O^{cp}})\cup\biguplus_{i\in I}C_i)\backslash\biguplus_{j\in J}C^*_j,$$ 
where $\text{BBY}_{\sigma,\sigma^*}^{-1}(\overrightarrow{O^{cp}})$ is a spanning tree. 
So $I=\emptyset$ $\Rightarrow$ $\varphi_{\sigma,\sigma^*}(\overrightarrow{O})$ is a forest. The inverse (and hence the converse) is true because $(\biguplus_{i\in I}C_i)\cap(\biguplus_{j\in J}C^*_j)=\emptyset$. 

The other equivalence is due to Corollary~\ref{cor0}(3).

(2) The proof is similar. 
\end{proof}

The next proposition shows that $\varphi_{\sigma,\sigma^*}$ is ``locally bijective'', which means if we fix a partial orientation or a subgraph, then the restricted $\varphi_{\sigma,\sigma^*}$ is still bijective. 

\begin{proposition}\label{local}

(1) Let $\overrightarrow{O_P}$ be a partial orientation of $G$, where $P$ is the underlying edges of $\overrightarrow{O_P}$. Then the map 

\begin{eqnarray*}
\varphi_{\sigma,\sigma^*,\overrightarrow{O_P}}:\{\text{orientations of }E\backslash P\} & \longrightarrow & \{\text{subsets of }E\backslash P\} \\
\overrightarrow{O_{E\backslash P}} & \mapsto & \varphi_{\sigma,\sigma^*}(\overrightarrow{O_{E\backslash P}}\cup\overrightarrow{O_P})\backslash P
\end{eqnarray*}
is a bijection. In other words, if one restricts $\varphi_{\sigma,\sigma^*}$ to the orientations containing $\overrightarrow{O_P}$ and ignores the edges in $P$, then one gets a bijection.

(2) Let $E_c$ and $E_d$ be two disjoint subsets of $E$. Then the map 

\begin{eqnarray*}
\varphi^{-1}_{\sigma,\sigma^*, E_c, E_d}: \{\text{subsets of }E\backslash (E_c\cup E_d)\} & \longrightarrow & \{\text{orientations of }E\backslash (E_c\cup E_d)\} \\
H & \mapsto & \varphi^{-1}_{\sigma,\sigma^*}(H\cup E_c)\text{ restricted to }E\backslash (E_c\cup E_d)
\end{eqnarray*}
is a bijection. In other words, if one restricts $\varphi^{-1}_{\sigma,\sigma^*}$ to the spanning subgraphs that include $E_c$ and exclude $E_d$, and ignores the edges in $E_c\cup E_d$, then one gets a bijection. 

\end{proposition}

\begin{proof}
By Proposition~\ref{propsubset}, these two maps are injective. Clearly for each map the domain and codomain have the same cardinality, so they are bijections.   
\end{proof}

\begin{example}
Here we illustrate the maps $\varphi_{\sigma,\sigma^*,\overrightarrow{O_P}}$ and $\varphi^{-1}_{\sigma,\sigma^*, E_c, E_d}$ of Proposition~\ref{local}. 
We choose the bijection $\varphi_{\sigma,\sigma^*}$ to be the one in Example~\ref{ex1}; see Figure~\ref{F5}. In Figure~\ref{F51}(a), we choose $\overrightarrow{O_P}$ to be one bottom arc directed from left to right. Then there are four orientations in Figure~\ref{F5} extending $\overrightarrow{O_P}$. By ignoring this bottom edge, we get four different subsets and four different orientations of $E\backslash P$. The naturally inherited correspondence between them is $\varphi_{\sigma,\sigma^*,\overrightarrow{O_P}}$. 
To fix a partial subgraph, we need to determine which edges are in it and which edges are not in it. In Figure~\ref{F51}(b), we include the left edge and exclude the right edge. Then there are two configurations satisfying this condition. By ignoring these two edges in the two configurations, we get the bijection $\varphi^{-1}_{\sigma,\sigma^*, E_c, E_d}$. 

\begin{figure}[h]
            \centering
            \includegraphics[scale=1]{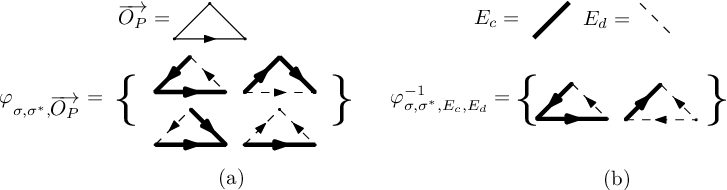}
            \caption{(a) We restrict the bijection $\varphi_{\sigma,\sigma^*}$ in Figure~\ref{F5} to the orientations extending $\protect\overrightarrow{O_P}$ and get $\varphi_{\sigma,\sigma^*,\protect\overrightarrow{O_P}}$. Technically speaking, the bottom arc in the four configurations of $\varphi_{\sigma,\sigma^*,\protect\overrightarrow{O_P}}$ should be erased. Here we still draw it to emphasize the relation between $\varphi_{\sigma,\sigma^*}$ and  $\varphi_{\sigma,\sigma^*,\protect\overrightarrow{O_P}}$. (b) We restrict the bijection $\varphi_{\sigma,\sigma^*}$ in Figure~\ref{F5} to the subgraphs that include the left edge and exclude the right edge and get $\varphi^{-1}_{\sigma,\sigma^*, E_c, E_d}$. }
            \label{F51}
\end{figure}

\end{example}

\section{The geometric proof of Theorem~\ref{th}(2)}\label{geometric}

This section aims at giving a geometric proof of Theorem~\ref{th}(2), which extends the geometric proof of Proposition~\ref{tree} in \cite{BBY}. We also explain how to prove Theorem~\ref{th}(3) by duality. 

The geometric approach in the section is independent of the combinatorial one in Section~\ref{combinatorial} and Section~\ref{tiling}.

There are many notations in this section. For the reader's convenience, we list them here. They will still be defined where they first appear except the last two. 

\noindent\textbf{Notation.} 

$\psi$: the projection map from the cube $[0,1]^E$ to the zonotope $Z_D$, which is a restriction of the linear map $D$. 

$\mu$: the section of $\psi$ whose image is the set of $\sigma$-compatible continuous orientations. 

$S_\sigma$:  the set of $\sigma$-compatible continuous orientations.

$\text{Face}_\sigma(T)$: the set of the continuous orientations where each edge $e\notin T$ is oriented according to $\sigma(C(T,e))$, which is a face of the cube $[0,1]^E$. 

$Z_\sigma(T):=\psi(\text{Face}_\sigma(T))$, which is a parallelepiped in the zonotope $Z_D$. 

$\widetilde{S_\sigma}:=S_\sigma+\im_\mathbb{Z}(D^T)$, which is a subset of $\mathbb{R}^E$ containing $S_\sigma$.

$\widetilde{\psi}$: a bijection from $\widetilde{S_\sigma}$ to $\mathbb{R}^r$, which is a restriction of the linear map $D$.  

$\text{hoc}_{\sigma,\sigma^*}(T)$:  the set of the continuous orientations where each edge $e\notin T$ is oriented according to $\sigma(C(T,e))$ and each edge $e\in T$ is oriented according to $\sigma(C(T,e))$ or bi-oriented, which is a half-open cell in $\text{Face}_\sigma(T)$.

$\widehat{Z_D}$: the $(q_1,\cdots,q_n)$-dilate of $Z_D$.

$\widehat{S_\sigma}$: the $(q_1,\cdots,q_n)$-dilate of $S_\sigma$.

$\widehat{\psi}$: a map from $\widehat{S_\sigma}$ to $\widehat{Z_D}$, which extends $\psi$.

$A+B:=\{a+b:a\in A, b\in B\}$, where $A$ and $B$ are two subsets of a linear space. 

$A\biguplus B$: the disjoint union of two sets $A$ and $B$. 

\subsection{An introduction to the geometric proof of Proposition~\ref{tree} in \cite{BBY}}
Fix a graph $G$, a reference orientation, an acyclic cycle signature $\sigma$, and an acyclic cocycle signature $\sigma^*$. 

We identity the set of \emph{continuous orientations} of $G$ with the cube $[0,1]^E$. In a continuous orientation, if the $e$-th coordinate is in $(0,1)$, then we say the edge $e$ is \emph{bi-oriented}.

We recall some definitions concerning the continuous orientations in \cite{BBY}.

A \emph{continuous cycle reversal} with respect to a directed cycle $\overrightarrow{C}$ replaces a continuous orientation $\overrightarrow{O}$ with $\overrightarrow{O}-\epsilon\overrightarrow{C}$ for some $\epsilon>0$ provided the new continuous orientation is still in $[0,1]^E$. Here $\overrightarrow{C}$ is not necessarily a directed cycle in $\overrightarrow{O}$. The equivalence relation generated by continuous cycle reversals defines the \emph{continuous cycle reversal (equivalence) classes}. 

A continuous orientation $\overrightarrow{O}$ is called \emph{$\sigma$-compatible} if one cannot flow $\overrightarrow{O}$ along any directed cycle in $\sigma$, i.e, $\overrightarrow{O}+\epsilon\sigma(C)\notin [0,1]^E$ for any $\epsilon >0$ and any cycle $C$. 

The \emph{continuous cocycle (equivalence) classes} and \emph{$\sigma^*$-compatible continuous orientations} are defined similarly.

The next lemma gives a geometric description of a cycle signature being acyclic. 
\begin{lemma}\cite[Lemma 3.1.1]{BBY}\label{signature}
Let $\sigma$ be a cycle signature of $G$. Then $\sigma$ is acyclic if and only if there exists $w\in\mathbb{R}^E$ such that $\langle w, \sigma(C)\rangle>0$ for each cycle $C$ of $G$.
\end{lemma}

Recall that $D$ is the modified incidence matrix of the reference orientation with one row removed. 

We restrict the linear map $D$ to the continuous orientations $[0,1]^E$. The image is the \emph{zonotope} $Z_D=\{\sum\limits_{i=1}^{|E|}c_iv_i:0\leq c_i \leq 1\}$, where $v_i$'s are the columns of $D$. We call this map $\psi$, i.e.,
\begin{eqnarray*}
\psi:[0,1]^E & \longrightarrow & Z_D \\
\overrightarrow{O} & \mapsto & D\cdot \overrightarrow{O},
\end{eqnarray*}
where $\overrightarrow{O}$ is viewed as a column vector.

The next proposition gives a simple description of the \emph{fibers} of $\psi$. 

\begin{proposition}\cite[Prop. 3.1.4 and Prop. 4.1.3]{BBY}\label{fiber}

(1) The map $\psi$ gives a bijection between continuous cycle reversal classes of continuous orientations of $G$ and points of the zonotope $Z_D$. In other words, two continuous orientations $\overrightarrow{O_1}$ and $\overrightarrow{O_2}$ are in the same continuous cycle reversal class if and only if $\overrightarrow{O_1}-\overrightarrow{O_2}\in\ker_\mathbb{R}(D)$.

(2) The map $\psi$ induces a bijection between (discrete) cycle reversal classes of $G$ and lattice points of the zonotope $Z_D$. 
\end{proposition}

The next result is the continuous version of Proposition~\ref{sigmacompatible1}.

\begin{proposition}\cite[Prop. 3.2.1]{BBY}\label{sigmacompatible2}
Let $\sigma$ be an acyclic cycle signature. Then each continuous cycle reversal class contains a unique $\sigma$-compatible continuous orientation.  
\end{proposition}

The similar results to Lemma~\ref{signature}, Proposition~\ref{fiber}, and Proposition~\ref{sigmacompatible2} hold for cocycles.

By \cite[Remark 3.2.2]{BBY}, the map
\begin{eqnarray*}
\mu:Z_D & \longrightarrow & [0,1]^E \\
z & \mapsto & \text{the unique }\sigma\text{-compatible continuous orientation in }\psi^{-1}(z)
\end{eqnarray*}
is a continuous \emph{section} to the map $\psi$, i.e., $\mu$ is continuous and $\psi\circ\mu$ is the identity map. 
By abusing language, we also call the image of a section a section. 

We denote $\im(\mu)$, the set of $\sigma$-compatible continuous orientations, by $S_\sigma$.

By Proposition~\ref{fiber}, Proposition~\ref{sigmacompatible1}, and Proposition~\ref{sigmacompatible2}, we have the following corollary. 

\begin{corollary}\label{section}
The map $\psi$ restricted to $S_\sigma$ is a bijection to $Z_D$. Moreover, $\psi$ restricted to the lattice points of $S_\sigma$, i.e., the $\sigma$-compatible orientations, is a bijection to the lattice points of $Z_D$. 
\end{corollary}

For each spanning tree $T$, define the face $\text{Face}_\sigma(T)$ of $[0,1]^E$ to be the set of the continuous orientations where each edge $e\notin T$ is oriented according to $\sigma(C(T,e))$. Define the parallelepiped $Z_\sigma(T)=\psi(\text{Face}_\sigma(T))\subseteq Z_D$. 

The following result shows a nice structure of $Z_D$.

\begin{proposition}\cite[Prop. 3.4.1]{BBY}\label{polyhedralzonotope}
Fix an acyclic cycle signature $\sigma$ of $G$. 

(1) $Z_D=\bigcup_{T\in\mathcal{T}(G)}Z_\sigma(T)$.

(2) For two different spanning trees $T_1$ and $T_2$, $Z_\sigma(T_1)$ and $Z_\sigma(T_2)$ are different and their intersection is either a common face or the empty set. 

(3) The vertices of $Z_\sigma(T)$'s correspond via $\mu$ to the $\sigma$-compatible (discrete) orientations (and hence by Corollary~\ref{section}, the vertices of $Z_\sigma(T)$'s are exactly the lattice points of $Z_D$). 
\end{proposition}

We also need the counterpart of Proposition~\ref{polyhedralzonotope}(1) for $S_\sigma$.

\begin{proposition}\cite{BBY}\label{verticesofS}
$S_\sigma=\bigcup_{T\in\mathcal{T}(G)}\text{Face}_\sigma(T)$. 
\end{proposition}

The fact that $\text{Face}_\sigma(T)\subseteq S_\sigma$ is proved implicitly in \cite[Prop. 3.3.1(2)]{BBY}; see also Lemma~\ref{facelemma}. The other direction $S_\sigma\subseteq\bigcup_{T\in\mathcal{T}(G)}\text{Face}_\sigma(T)$ is proved implicitly in \cite[Prop. 3.4.1]{BBY}.

\begin{figure}[h]
            \centering
            \includegraphics[scale=1]{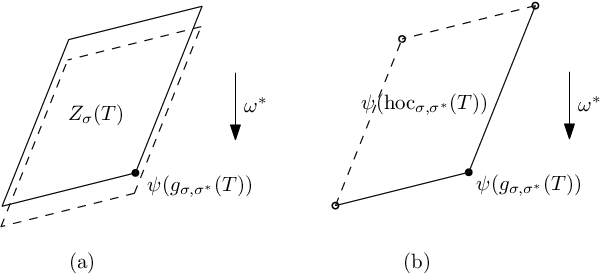}
            \caption{(a): We shift the parallelepiped $Z_\sigma(T)$ along $w^*$, then it will cover the vertex $\psi(\text{BBY}_{\sigma,\sigma^*}(T))$. (b): The half-open cell $\text{hoc}_{\sigma,\sigma^*}(T)$ defined in Subsection~\ref{hoc}. }
            \label{F0}
\end{figure}

Now we consider shifting the parallelepiped $Z_\sigma(T)$ in the ambient Euclidean space of $Z_D$ along a generic vector $w^*$; see Figure~\ref{F0}(a). For sufficiently small positive $\epsilon$, the image of $Z_\sigma(T)$ under the shifting map $v\mapsto v+\epsilon w^*$ contains a unique vertex of $Z_\sigma(T)$ corresponding via $\mu$ to a $\sigma$-compatible orientation. This defines a map that sends a spanning tree to a $\sigma$-compatible orientation. Fix an acyclic cocycle signature $\sigma^*$. By \cite[Theorem 3.5.3]{BBY}, there exists\footnote{The paper \cite{BBY} works with the row zonotope instead of the zonotope when explaining the geometric interpretation of $\text{BBY}_{\sigma,\sigma^*}$. By \cite[Lemma 3.5.1]{BBY}, the two zonotopes differ by the linear transformation $D$, so the geometric idea work for both zonotopes. The vector $w^*$ can be chosen to be $D\cdot w'$, where $w'$ is the shifting vector used for the row zonotope in \cite[Lemma 3.5.2]{BBY}. For simplicity, we only claim the existence of such a vector.} a vector $w^*$ (depending on $\sigma^*$) such that the map is $\text{BBY}_{\sigma,\sigma^*}$ defined in Proposition~\ref{tree}. By \cite[Theorem 3.5.5]{BBY}, $\text{BBY}_{\sigma,\sigma^*}$ is a bijection to the set of  $(\sigma,\sigma^*)$-compatible orientations, which is Proposition~\ref{tree}. We summarize this paragraph with the following proposition. 

\begin{proposition}\cite{BBY}\label{bbymain}
Fix acyclic signatures $\sigma$ and $\sigma^*$ of $G$. Then there exists a vector $w^*$ such that, for sufficiently small positive $\epsilon$ and any spanning tree $T$, the image of $Z_\sigma(T)$ under the shifting map $v\mapsto v+\epsilon w^*$ contains a unique vertex of $Z_\sigma(T)$ corresponding (via $\mu$) to the $(\sigma,\sigma^*)$-compatible orientation $\text{BBY}_{\sigma,\sigma^*}(T)$ defined in Proposition~\ref{tree}. Moreover, the map $\text{BBY}_{\sigma,\sigma^*}$ is bijective from $\mathcal{T}(G)$ to $\{(\sigma,\sigma^*)\text{-compatible orientations}\}$.
\end{proposition}

See Section~\ref{intro} for an example.

\subsection{An extension of the section $S_\sigma$}
The target of this subsection is to prove that $$\widetilde{S_\sigma}:=S_\sigma+\im_\mathbb{Z}(D^T)$$ is the image of a section to the linear map $D$ that extends $\mu$. Intuitively, one may imagine $\widetilde{S_\sigma}$ as Figure~\ref{F3}(a), where the copies of $Z_D$ are viewed as copies of $S_\sigma$.

\begin{Rem}
As introduced in Section~\ref{intro}, we want to make use of the fact that copies of the zonotope $Z_D$ tile the ambient space $\im_\mathbb{R}(D)$. For technical reasons, we make the tiling for the section $\widetilde{S_\sigma}$ and transfer the tiling to $\im_\mathbb{R}(D)$ via $\psi$. So we do not make use of the result \cite[Theorem 1]{DG} for the space tiling zonotopes directly. 
\end{Rem}

We start with a lemma, where intuitively $[0,1]^E+\im_\mathbb{Z}(D^T)$ looks like Figure~\ref{F3}(a) if one imagine copies of $Z_D$ as copies of $[0,1]^E$. 

\begin{lemma}\label{cubepile}
$\mathbb{R}^E=[0,1]^E+\im_\mathbb{Z}(D^T)+\ker_\mathbb{R}(D)$. 
\end{lemma}
\begin{proof}
Let $\overrightarrow{O}\in\mathbb{R}^E$. Denote the coordinate of $\overrightarrow{O}$ corresponding to an edge $e$ by $e(\overrightarrow{O})$. Define the distance function $N(\overrightarrow{O})=\max_{e\in E}|e(\overrightarrow{O})-1/2|$. Note that $\overrightarrow{O}\in [0,1]^E$ if and only if $N(\overrightarrow{O})\leq 1/2$.
This function measures the ``distance'' between $\overrightarrow{O}$ and the cube $[0,1]^E$. 

We first show that (i) if $\overrightarrow{O}\notin [0,1]^E$, then $N(\overrightarrow{O}+\overrightarrow{v^*}+\overrightarrow{v})<N(\overrightarrow{O})$ for some $\overrightarrow{v^*}\in\im_\mathbb{Z}(D^T)$ and some $\overrightarrow{v}\in\ker_\mathbb{R}(D)$. Then we show that (ii) when a vector $\overrightarrow{O'}$ varies in 
$\overrightarrow{O}+\im_\mathbb{Z}(D^T)+\ker_\mathbb{R}(D)$,  $N(\overrightarrow{O'})$ can attain its minimum value. Hence the lemma holds.

To prove (i), we assume $\overrightarrow{O}\notin [0,1]^E$. Let $P$ be the set of edges $e$ such that $|e(\overrightarrow{O})-1/2|=N(\overrightarrow{O})$. Note that $e(\overrightarrow{O})>1$ or $e(\overrightarrow{O})<0$. We construct a partial orientation $\overrightarrow{P}$ as follows. If $e\in P$ satisfies $e(\overrightarrow{O})>1$, then we put the arc $\overrightarrow{e}$ oriented by the reference orientation into $\overrightarrow{P}$; if $e\in P$ satisfies $e(\overrightarrow{O})<0$, then we put the arc $\overrightarrow{e}$ oriented against the reference orientation into $\overrightarrow{P}$. (Recall that the coordinate of the reference orientation is $(1,1,\cdots,1)$.)

Now we apply Lemma~\ref{3-painting} to the data $\overrightarrow{P}$, $E_c=E\backslash P$, $E_d=\emptyset$, and any arc $\overrightarrow{e}\in\overrightarrow{P}$. If $\overrightarrow{e}$ belongs a directed cocycle $\overrightarrow{C^*}\in\im_\mathbb{Z}(D^T)$, then the new vector $\overrightarrow{O'}=\overrightarrow{O}-\overrightarrow{C^*}$ satisfies $N(\overrightarrow{O'})\leq N(\overrightarrow{O})$ and the number of edges $e'$ such that $|e'(\overrightarrow{O'})-1/2|=N(\overrightarrow{O})$ decreases. In the other case that $\overrightarrow{e}$ belongs a directed cycle $\overrightarrow{C}\in\ker_\mathbb{R}(D)$, then for a sufficiently small positive number $\epsilon$, the new vector $\overrightarrow{O'}=\overrightarrow{O}-\epsilon\overrightarrow{C}$ satisfies $N(\overrightarrow{O'})\leq N(\overrightarrow{O})$ and the number of edges $e'$ such that $|e'(\overrightarrow{O'})-1/2|=N(\overrightarrow{O})$ decreases. So we can apply the same process until $N(\overrightarrow{O'})<N(\overrightarrow{O})$ and hence (i) holds. 

It remains to prove (ii). Note that the function $N$ is continuous, so we will make use of compactness. Consider the orthogonal decomposition $\ker_\mathbb{R}(D)\oplus\im_\mathbb{R}(D^T)$ of $\mathbb{R}^E$ and its topology. Note that $\im_\mathbb{Z}(D^T)$ is closed in $\im_\mathbb{R}(D^T)$ and $\ker_\mathbb{R}(D)$ is closed in $\ker_\mathbb{R}(D)$, so $\overrightarrow{O}+\im_\mathbb{Z}(D^T)+\ker_\mathbb{R}(D)$ is closed in $\mathbb{R}^E$. Let $r=N(\overrightarrow{O})$. Then the set
$(\overrightarrow{O}+\im_\mathbb{Z}(D^T)+\ker_\mathbb{R}(D))\cap \{\overrightarrow{O'}\in\mathbb{R}^E:N(\overrightarrow{O'})\leq r\}$ is compact and hence (ii) holds.

\end{proof}

\begin{corollary}\label{decomposition2}
$\mathbb{R}^E=S_\sigma+\im_\mathbb{Z}(D^T)+\ker_\mathbb{R}(D)$.
\end{corollary}
\begin{proof}
By Lemma~\ref{cubepile}, any vector in $\mathbb{R}^E$ can be written as $\overrightarrow{O}+\overrightarrow{v^*}+\overrightarrow{v}$, where $\overrightarrow{O}\in [0,1]^E$, $\overrightarrow{v^*}\in\im_\mathbb{Z}(D^T)$, and $\overrightarrow{v}\in\ker_\mathbb{R}(D)$. 
By Proposition~\ref{sigmacompatible2}, $\overrightarrow{O}$ and some $\overrightarrow{O'}\in S_\sigma$ are in the same continuous cycle reversal class. By Proposition~\ref{fiber}(1), $\overrightarrow{O}-\overrightarrow{O'}\in\ker_\mathbb{R}(D)$. So we have the desired sum  $\overrightarrow{O'}+\overrightarrow{v^*}+(\overrightarrow{v}+\overrightarrow{O}-\overrightarrow{O'})$.
\end{proof}

\begin{proposition}\label{decomposition3}
$\mathbb{R}^E=\widetilde{S_\sigma}+\ker_\mathbb{R}(D)$. Moreover, the sum is \emph{direct} in the sense that any vector in $\mathbb{R}^E$ can be written uniquely as the sum of a vector in $\widetilde{S_\sigma}$ and a vector in $\ker_\mathbb{R}(D)$. 
\end{proposition}
\begin{proof}

Recall that 
$\widetilde{S_\sigma}=S_\sigma+\im_\mathbb{Z}(D^T)$. By Corollary~\ref{decomposition2}, we get $\mathbb{R}^E=\widetilde{S_\sigma}+\ker_\mathbb{R}(D)$. 

To prove uniqueness, we let $$(\overrightarrow{O_1}+\overrightarrow{C_1^*})+\overrightarrow{C_1}=(\overrightarrow{O_2}+\overrightarrow{C_2^*})+\overrightarrow{C_2},$$ where  $\overrightarrow{O_k}\in S_\sigma\subseteq [0,1]^E$, $\overrightarrow{C_k^*}\in\im_\mathbb{Z}(D^T)$, and $\overrightarrow{C_k}\in\ker_\mathbb{R}(D)$ for $k=1,2$. It suffices to show that $\overrightarrow{O_1}+\overrightarrow{C_1^*}=\overrightarrow{O_2}+\overrightarrow{C_2^*}$. Denote $\overrightarrow{C^*}=\overrightarrow{C_2^*}-\overrightarrow{C_1^*}$ and $\overrightarrow{C}=\overrightarrow{C_2}-\overrightarrow{C_1}$. Then we have 
 $$\overrightarrow{O_1}-\overrightarrow{O_2}-\overrightarrow{C^*}=\overrightarrow{C},$$ where $\overrightarrow{C^*}\in\im_\mathbb{Z}(D^T)$ and $\overrightarrow{C}\in\ker_\mathbb{R}(D)$. We want to show $\overrightarrow{O_1}-\overrightarrow{C^*}=\overrightarrow{O_2}$.

First we show that $\overrightarrow{C^*}$ is a sum of disjoint directed cocycles, and the discretely oriented edges in $\overrightarrow{O_1}$ contains the disjoint cocycles $\overrightarrow{C^*}$. Because $\ker_\mathbb{R}(D)\perp \im_\mathbb{R}(D^T)$, we get $\langle\overrightarrow{O_1}-\overrightarrow{O_2}-\overrightarrow{C^*}, \overrightarrow{C^*}\rangle=0$ and hence 
$$\langle\overrightarrow{O_1}-\overrightarrow{O_2}, \overrightarrow{C^*}\rangle=\langle\overrightarrow{C^*}, \overrightarrow{C^*}\rangle.$$ Observe that $\overrightarrow{O_1}-\overrightarrow{O_2}\in[-1,1]^E$ and  $\overrightarrow{C^*}\in\mathbb{Z}^E$ imply $\langle\overrightarrow{O_1}-\overrightarrow{O_2}, \overrightarrow{C^*}\rangle\leq\langle\overrightarrow{C^*}, \overrightarrow{C^*}\rangle$, where the equality holds if and only if $\overrightarrow{C^*}$ is a $\{0,\pm 1\}$-vector and the coordinates of $\overrightarrow{O_1}-\overrightarrow{O_2}$ coincide with  $\overrightarrow{C^*}$ on the support of $\overrightarrow{C^*}$ . So, $\overrightarrow{C^*}$ is a sum of disjoint directed cocycles (by Lemma~\ref{conformal}), and  $\overrightarrow{O_1}$ contains $\overrightarrow{C^*}$. (By ``$\overrightarrow{O_1}$ contains $\overrightarrow{C^*}$'', we mean that when restricted to the underlying edges of $\overrightarrow{C^*}$, the continuous orientation $\overrightarrow{O_1}$ is discrete and coincides with $\overrightarrow{C^*}$. However, their coordinates are not equal in general, because when an arc disagrees with the reference orientation, its coordinate is $0$ in $\overrightarrow{O}$ but $-1$ in $\overrightarrow{C^*}$.)

Let $\overrightarrow{O_3}:=\overrightarrow{O_1}-\overrightarrow{C^*}$, which is the continuous orientation obtained by reversing $\overrightarrow{C^*}$ in $\overrightarrow{O_1}$. Because $\overrightarrow{O_1}\in S_\sigma$,  we have $\overrightarrow{O_3}\in S_\sigma$. By $\overrightarrow{O_3}-\overrightarrow{O_2}=\overrightarrow{C}\in\ker_\mathbb{R}(D)$ and Proposition~\ref{fiber}(1), $\overrightarrow{O_2}$ and $\overrightarrow{O_3}$ are in the same continuous cycle reversal class. Then by Proposition~\ref{sigmacompatible2}, $\overrightarrow{O_3}=\overrightarrow{O_2}$. 

\end{proof}

\begin{corollary}\label{tinycor}
$\widetilde{S_\sigma}\cap [0,1]^E=S_\sigma$.
\end{corollary}
\begin{proof}
$\widetilde{S_\sigma}\cap [0,1]^E\supseteq S_\sigma$ is trivial.

For the other direction, let $\overrightarrow{O}\in \widetilde{S_\sigma}$ be a continuous orientation. By Proposition~\ref{sigmacompatible2} and Proposition~\ref{fiber}(1), $\overrightarrow{O}=\overrightarrow{O}'+\overrightarrow{C}$ for some $\overrightarrow{O}'\in S_\sigma\subseteq\widetilde{S_\sigma}$ and $\overrightarrow{C}\in\ker_\mathbb{R}(D)$. By Proposition~\ref{decomposition3}, we compare both sides of $\overrightarrow{O}+0=\overrightarrow{O}'+\overrightarrow{C}$ and get $\overrightarrow{O}=\overrightarrow{O}'\in S_\sigma$.

\end{proof}

By Proposition~\ref{decomposition3}, we get a section to the linear map $D$ which sends a point $p$ in $\im_\mathbb{R}(D)$ to the unique point in the intersection of its preimage and $\widetilde{S_\sigma}$. It is easy to see that this section extends $\mu$. It seems more convenient to work with a bijection than a section. Now we define it. We denote the rank of $D$ by $r$, and then  $\im_\mathbb{R}(D)=\mathbb{R}^r\supseteq Z_D$.

\begin{proposition}\label{bigsection}
The map  
\begin{eqnarray*}
\widetilde{\psi}:\widetilde{S_\sigma} & \longrightarrow & \mathbb{R}^r \\
\overrightarrow{O} & \mapsto & D\cdot \overrightarrow{O}
\end{eqnarray*}
is a bijection, where $\overrightarrow{O}$ is viewed as a column vector. Moreover, $\widetilde{\psi}$ sends $S_\sigma$ to $Z_D$ and preserve lattice points, i.e., $\overrightarrow{O}\in\mathbb{Z}^E$ if and only if $\widetilde{\psi}(\overrightarrow{O})\in\mathbb{Z}^r$.
\end{proposition}
\begin{proof}
By Proposition~\ref{decomposition3}, $\widetilde{\psi}$ is a bijection. By Corollary~\ref{section}, $\widetilde{\psi}$ sends $S_\sigma$ to $Z_D$. It is trivial that $\overrightarrow{O}\in\mathbb{Z}^E$ implies $\widetilde{\psi}(\overrightarrow{O})\in\mathbb{Z}^r$. 
For the other direction, write $\overrightarrow{O}=\overrightarrow{O'}+\overrightarrow{C^*}$, 
where $\overrightarrow{O'}\in S_\sigma$ and $\overrightarrow{C^*}\in\im_\mathbb{Z}(D^T)$. Then $\widetilde{\psi}(\overrightarrow{O})=\psi(\overrightarrow{O'})+D\cdot\overrightarrow{C^*}$ and hence $\psi(\overrightarrow{O'})\in\mathbb{Z}^r$. By Corollary~\ref{section}, $\overrightarrow{O'}\in\mathbb{Z}^E$, so $\overrightarrow{O}\in\mathbb{Z}^E$. 
\end{proof}

\subsection{Polyhedral subdivision and half-open decomposition of $\widetilde{S_\sigma}$}\label{hoc}
By a \emph{polyhedral subdivision} of a subset $S\subseteq \mathbb{R}^n$ we mean writing $S$ as the union of a collection of $n$-dimensional polytopes where any two of the polytopes intersect in their common face. By Proposition~\ref{polyhedralzonotope}, we have a polyhedral subdivision of $Z_D$ (and hence $S_\sigma$ via $\mu$). Now we generalize it to $\widetilde{S_\sigma}$.

\begin{proposition}\label{subdivision1}
Fix an acyclic cycle signature $\sigma$ of $G$, then $$\widetilde{S_\sigma}=\bigcup_{T\in\mathcal{T}(G),\overrightarrow{C^*}\in\im_\mathbb{Z}(D^T) }(\text{Face}_\sigma(T)+\overrightarrow{C^*}).$$
Moreover, for any two different pairs $(T_1, \overrightarrow{C^*_1})$ and $(T_2, \overrightarrow{C^*_2})$ in $\mathcal{T}(G)\times\im_\mathbb{Z}(D^T)$, $\text{Face}_\sigma(T_1)+\overrightarrow{C^*_1}$ and $\text{Face}_\sigma(T_2)+\overrightarrow{C^*_2}$ are different and their intersection is either a common face or the empty set.

\end{proposition}
\begin{proof}
By Proposition~\ref{verticesofS}, we have  $S_\sigma=\bigcup_{T\in\mathcal{T}(G)}\text{Face}_\sigma(T)$. Because $\widetilde{S_\sigma}=S_\sigma+\im_\mathbb{Z}(D^T)$, the equality in the proposition holds. 

It remains to show that this is a polyhedral subdivision. Without loss of generality, assume $\overrightarrow{C^*_1}=0$. Let $\overrightarrow{O}\in\text{Face}_\sigma(T_1)\cap(\text{Face}_\sigma(T_2)+\overrightarrow{C^*_2})$. Because $\overrightarrow{O},\overrightarrow{O}-\overrightarrow{C^*_2}\in [0,1]^E$, $\overrightarrow{C^*_2}$ is a $\{0,\pm 1\}$-vector, which is a sum of disjoint directed cocycles by Lemma~\ref{conformal}, and the continuous orientation $\overrightarrow{O}$ contains these directed cocycles. 

So a continuous orientation $\overrightarrow{O}$ is in $\text{Face}_\sigma(T_1)\cap(\text{Face}_\sigma(T_2)+\overrightarrow{C^*_2})$ if and only if (i) $\overrightarrow{O}$ is in $\text{Face}_\sigma(T_1)$, (ii) $\overrightarrow{O}$ contains $\overrightarrow{C^*_2}$, and (iii) after reversing $\overrightarrow{C^*_2}$ in $\overrightarrow{O}$, any edge $e\notin T_2$ is oriented according to $\sigma(C(T_2,e))$. Note that a face of the parallelepiped $\text{Face}_\sigma(T_1)$ is a set of continuous orientations $\overrightarrow{O}$ in $\text{Face}_\sigma(T_1)$ where the coordinates corresponding to some edges in $T_1$ are $0$ or $1$. The conditions (ii) and (iii) also mean some coordinates in $\overrightarrow{O}$ are $0$ or $1$. So, if there is any contradiction among (i), (ii), and (iii), then the intersection is the empty set; otherwise the intersection is a face of $\text{Face}_\sigma(T_1)$. Moreover, the face cannot be $\text{Face}_\sigma(T_1)$ because any cocycle intersects any spanning tree. Similarly, we can prove $\text{Face}_\sigma(T_1)\cap(\text{Face}_\sigma(T_2)+\overrightarrow{C^*_2})$ is the empty set or a proper face of $\text{Face}_\sigma(T_2)+\overrightarrow{C^*_2}$.

\end{proof}

\begin{corollary}\label{subdivision2}
Fix an acyclic cycle signature $\sigma$ of $G$, then $$\mathbb{R}^r=\bigcup_{T\in\mathcal{T}(G),\overrightarrow{C^*}\in\im_\mathbb{Z}(D^T) }(Z_\sigma(T)+D\cdot\overrightarrow{C^*}).$$
Moreover, for any two different pairs $(T_1, \overrightarrow{C^*_1})$ and $(T_2, \overrightarrow{C^*_2})$ in $\mathcal{T}(G)\times\im_\mathbb{Z}(D^T)$, $Z_\sigma(T_1)+D\cdot\overrightarrow{C^*_1}$ and $Z_\sigma(T_2)+D\cdot\overrightarrow{C^*_2}$ are different and their intersection is either a common face or the empty set.

\end{corollary}
\begin{proof}
Note that $\psi$ sends the $r$ linearly independent generating edges of the parallelepiped $\text{Face}_\sigma(T)$ to the $r$ linearly independent generating edges of the parallelepiped $Z_\sigma(T)$, so $\psi$ sends a face of $\text{Face}_\sigma(T)$ to a face of $Z_\sigma(T)$. Hence this corollary is a direct consequence of Proposition~\ref{bigsection} and Proposition~\ref{subdivision1}.
\end{proof}

Remark that although Proposition~\ref{subdivision1} and Corollary~\ref{subdivision2} generalize Proposition~\ref{polyhedralzonotope}(2), we do not use it in the proofs. 

By a \emph{half-open decomposition} of a vector set $S$ we mean writing $S$ as a disjoint union of half-open cells. We now look for a half-open decomposition of $\mathbb{R}^r$ and hence of $\widetilde{S_\sigma}$. Fix $\sigma$ and $\sigma^*$, for each spanning tree $T$, define the half-open cell $\text{hoc}_{\sigma,\sigma^*}(T)$ to be the set of the continuous orientations where each edge $e\notin T$ is oriented according to $\sigma(C(T,e))$ and each edge $e\in T$ is oriented according to $\sigma(C(T,e))$ or bi-oriented. Recall in Proposition~\ref{bbymain} that there exists $w^*$ such that if we shift the parallelepiped $Z_\sigma(T)$ along $w^*$ a bit, then the unique vertex it will cover is $\psi(\text{BBY}_{\sigma,\sigma^*}(T))$. The half-open cell $\psi(\text{hoc}_{\sigma,\sigma^*}(T))$ can be obtained from $Z_\sigma(T)$ by removing all the closed facets not containing  $\psi(\text{BBY}_{\sigma,\sigma^*}(T))$; see Figure~\ref{F0}. By basic geometry, we have the following lemma.

\begin{lemma}\label{hoclemma}
The half-open cell $\psi(\text{hoc}_{\sigma,\sigma^*}(T))$ is the set points $p\in\mathbb{R}^r$ such that for sufficiently small positive $\epsilon$, $p-\epsilon w^*$ is in the interior of the parallelepiped $Z_\sigma(T)$. 
\end{lemma}

\begin{proposition}\label{hocprop}
Fix $\sigma$ and $\sigma^*$ of $G$. Then

(1) $\mathbb{R}^r=\biguplus_{T\in\mathcal{T}(G),\overrightarrow{C^*}\in\im_\mathbb{Z}(D^T)}(\psi(\text{hoc}_{\sigma,\sigma^*}(T))+D\cdot\overrightarrow{C^*})$; 

(2) $\widetilde{S_\sigma}=\biguplus_{T\in\mathcal{T}(G),\overrightarrow{C^*}\in\im_\mathbb{Z}(D^T) }(\text{hoc}_{\sigma,\sigma^*}(T)+\overrightarrow{C^*})$.

\end{proposition}

\begin{proof}
(1) Let $w^*$ be as in Proposition~\ref{bbymain} and $p\in\mathbb{R}^r$. By Corollary~\ref{subdivision2}, if we shift $p$ along $-w^*$ by a sufficiently small distance, then it will be in the interior of a unique parallelepiped $Z_\sigma(T)+D\cdot\overrightarrow{C^*}$. By Lemma~\ref{hoclemma}, $p\in\psi(\text{hoc}_{\sigma,\sigma^*}(T))+D\cdot\overrightarrow{C^*}$.

(2) Apply $\widetilde{\psi}^{-1}$ to the first identity, then we get the second one.  
\end{proof}

\subsection{Half-open decomposition of $S_\sigma$ and the induced map $\varphi_{\sigma,\sigma^*}$}\label{varphi}

In this subsection we will restrict the half-open decomposition of $\widetilde{S_\sigma}$ to $S_\sigma$ and give a geometric definition of the map $\varphi_{\sigma,\sigma^*}$ in Theorem~\ref{th}(2). Then we show that this map is the same as the one in Theorem~\ref{th}(2), which is defined in a combinatorial way. 

We start with some definitions. All the half-open cells we talk about in this subsection are in $\mathbb{R}^E$ and of the form $\prod_{e\in E} I_e$, where each $I_e$ is $\{a_e\}$, $(a_e-1,a_e]$, or $[a_e,a_e+1)$ for some integer $a_e$. We call them \emph{standard half-open cells} in $\mathbb{R}^E$. Any standard half-open cell contains a unique lattice point $(a_e)_{e\in E}$. We call the set $\{e\in E:I_e=(a_e-1,a_e]\text{ or }[a_e,a_e+1))\}$ the \emph{generating set} of the standard half-open cell. 

The following lemma is trivial. In order to state the results of this subsection in a better way, we change the sign here. 
\begin{lemma}\label{standard}
Fix $\sigma$ and $\sigma^*$ of $G$. Let $T\in\mathcal{T}(G)$ and $\overrightarrow{C^*}\in\im_\mathbb{Z}(D^T)$. 

(1) $\text{hoc}_{\sigma,\sigma^*}(T)-\overrightarrow{C^*}$ is a standard half-open cell. Its generating set is $T$ and it contains the unique lattice point $\text{BBY}_{\sigma,\sigma^*}(T)-\overrightarrow{C^*}$, where $\text{BBY}_{\sigma,\sigma^*}$ is the bijection defined in Proposition~\ref{tree}.

(2) If $(\text{hoc}_{\sigma,\sigma^*}(T)-\overrightarrow{C^*})\cap [0,1]^E\neq\emptyset$, then it is a standard half-open cell. Its generating set is a forest $F\subseteq T$ and it contains the unique lattice point $\text{BBY}_{\sigma,\sigma^*}(T)-\overrightarrow{C^*}$.
\end{lemma}

Now we get a half-open decomposition of $S_\sigma$ as follows. 

\begin{proposition}\label{hocprop2}
Fix $\sigma$ and $\sigma^*$ of $G$. Then
$$S_\sigma=\biguplus_{T\in\mathcal{T}(G),\overrightarrow{C^*}\in\im_\mathbb{Z}(D^T)}(\text{hoc}_{\sigma,\sigma^*}(T)-\overrightarrow{C^*})\cap [0,1]^E,$$ where each summand is a standard half-open cell if it is non-empty. 
\end{proposition}

\begin{proof}
This follows immediately from Proposition~\ref{hocprop}, Corollary~\ref{tinycor}, and Lemma~\ref{standard}. 

\end{proof}

Then we define the map $\varphi_{\sigma,\sigma^*}$ from a geometric point of view. For any $\sigma$-compatible (discrete) orientation $\overrightarrow{O}$, it is a lattice point of $S_\sigma$ and hence by Proposition~\ref{hocprop2}, there exists a unique standard half-open cell $(\text{hoc}_{\sigma,\sigma^*}(T)-\overrightarrow{C^*})\cap [0,1]^E$ containing it. Define $\varphi_{\sigma,\sigma^*}(\overrightarrow{O})$ to be the generating set of the half-open cell, which is a forest $F\subseteq T$. 

The next lemma characterizes the non-empty half-open cells in the decomposition of $S_\sigma$. 

\begin{lemma}\label{forestlemma}
Fix $\sigma$ and $\sigma^*$ of $G$. Let $T\in\mathcal{T}(G)$ and $\overrightarrow{C^*}\in\im_\mathbb{Z}(D^T)$. Then the following are equivalent:

(a) $(\text{hoc}_{\sigma,\sigma^*}(T)-\overrightarrow{C^*})\cap [0,1]^E\neq\emptyset$. 

(b) $\text{BBY}_{\sigma,\sigma^*}(T)-\overrightarrow{C^*}\in [0,1]^E$.

(c) $\overrightarrow{C^*}$ is a sum of disjoint directed cocycles and each directed cocycle is in the orientation $\text{BBY}_{\sigma,\sigma^*}(T)$.
\end{lemma}

\begin{proof}
Due to Lemma~\ref{standard}(1), (b) implies (a). 

Due to Lemma~\ref{standard}(2), (a) implies (b). 

It is trivial that (c) implies (b). 

To prove that (b) implies (c), we write (b) in terms of coordinates. Then we find that $\overrightarrow{C^*}$ has to be a $\{0,\pm 1\}$-vector and $\text{BBY}_{\sigma,\sigma^*}(T)$ restricted to the underlying edges of $\overrightarrow{C^*}$ is $\overrightarrow{C^*}$. By Lemma~\ref{conformal}, (c) holds. 

\end{proof}

\begin{corollary}\label{cor1}
If the standard half-open cell 
$H=(\text{hoc}_{\sigma,\sigma^*}(T)-\overrightarrow{C^*})\cap [0,1]^E$ in Proposition~\ref{hocprop2} is non-empty, then 

(1) $\overrightarrow{C^*}$ is a sum of disjoint directed cocycles $\{\overrightarrow{C^*_j}\}_{j\in J}$, where each $\overrightarrow{C^*_j}$ is in the orientation $\text{BBY}_{\sigma,\sigma^*}(T)$;

(2) the unique lattice point contained in $H$ is
$\text{BBY}_{\sigma,\sigma^*}(T)-\overrightarrow{C^*}$, which is an orientation obtained by reversing the directed cocycle $\overrightarrow{C_j^*}$'s in $\text{BBY}_{\sigma,\sigma^*}(T)$;

(3) and the generating set of $H$ is the forest $T\backslash (\biguplus_{j\in J}C_j^*)$.

\end{corollary}

\begin{proof}
The only part we have not proved is (3). 
Write the standard half-open cell $\text{hoc}_{\sigma,\sigma^*}(T)=\prod_{e\in E} I_e$, where each $I_e$ is $\{a_e\}$, $(a_e-1,a_e]$, or $[a_e,a_e+1)$ for some integer $a_e$. The generating set of $\text{hoc}_{\sigma,\sigma^*}(T)$ is $T$, consisting of edges $e$ such that $I_e$ is a half-open interval. To get $H$, we shift the cell by $-\overrightarrow{C^*}$ and restrict it to $[0,1]^E$. So, during this process a half-open interval $I_e$ becomes a point if and only if $I_e$ is shifted (and hence its interior part is moved out of the cube), which means $e\in\biguplus_{j\in J}C_j^*$. These edges are removed from the generating set, so (3) holds. 

\end{proof}

Now we claim that the map $\varphi_{\sigma,\sigma^*}$ is exactly the one in Theorem~\ref{th}(2). 
\begin{proposition}\label{forestgeometry}
Fix $\sigma$ and $\sigma^*$ of $G$. Let 
\begin{eqnarray*}
\varphi_{\sigma,\sigma^*}:\{\sigma\text{-compatible orientations}\}\longrightarrow \{\text{spanning forests}\}
\end{eqnarray*}
be the map that sends a $\sigma$-compatible orientation $\overrightarrow{O}$ to the generating forest of the unique standard half-open cell in Proposition~\ref{hocprop2} containing $\overrightarrow{O}$. Let $\overrightarrow{O^{cp}}$ be the ($\sigma,\sigma^*$)-compatible orientation in the cocycle reversal class containing $\overrightarrow{O}$ and $\{\overrightarrow{C_j^*}\}_{j\in J}$ be the disjoint directed cocycles in $\overrightarrow{O^{cp}}$ by reversing which we get $\overrightarrow{O}$. Then $\varphi_{\sigma,\sigma^*}(\overrightarrow{O})$ is the forest $\text{BBY}_{\sigma,\sigma^*}^{-1}(\overrightarrow{O^{cp}})\backslash (\biguplus_{j\in J}C_j^*)$. 
\end{proposition}

\begin{proof}
This is a direct consequence of Corollary~\ref{cor1} and Corollary~\ref{cor0}. 
\end{proof}

\subsection{Bijectivity of $\varphi_{\sigma,\sigma^*}$ and the generalized Ehrhart polynomial of a zonotope}\label{subsectionduality}

In this subsection we generalize the \emph{Ehrhart polynomial} of a zonotope to prove the bijectivity of the map $\varphi_{\sigma,\sigma^*}$ in Proposition~\ref{forestgeometry} (which is also the one in Theorem~\ref{th}(2)). 

We first recall the Ehrhart polynomial. For any positive integer $q$ and any convex polytope $\mathcal{P}\subseteq\mathbb{R}^r$ with integer vertices, denote by $E_\mathcal{P}(q)$ the number of lattice points in $q\cdot\mathcal{P}:=\{qp:p\in\mathcal{P}\}$ (including the boundary). 
Then $E_\mathcal{P}(q)$ is a polynomial function of $q$, called the \emph{Ehrhart polynomial} of $\mathcal{P}$. 

Let $A$ be a $r\times n$ integer matrix. Denote the columns of $A$ by $v_1, v_2, \cdots, v_n$ and the zonotope $\{\sum\limits_{i=1}^{n}c_iv_i:0\leq c_i \leq 1\}$ by $Z_A$. We rewrite $E_{Z_A}(q)$ as $E_A(q)$. 
Now we recall the following basic fact. 

\begin{lemma}\cite[Theorem 2.2]{S}\label{Ehrhart1}
The Ehrhart polynomial of $Z_A$ is
$$E_A(q)=\sum_Xh(X)q^{|X|},$$
where $X$ ranges over all linearly independent subsets of $\{v_1, v_2, \cdots, v_n\}$, and where $h(X)$ denotes the greatest common divisor of all minors of size $|X|$ of the matrix whose columns are the elements of $X$. 
\end{lemma}

From now on we take $A$ to be a \emph{totally unimodular matrix}, meaning every square submatrix has determinant $0$, $1$, or $-1$. We further assume that $\text{rank}(A)=r$. For example, the matrix $D$ we have been using in this section satisfies this condition \cite[Proposition 5.1.3]{O}. Then $h(X)=1$ for all linearly independent set $X$ and hence $E_A(q)=\sum_X q^{|X|}$. Note that when $A$ is taken to be $D$, the edges corresponding to $X$, which can also be viewed as the subscripts of the elements in $X$, form a forest. This is because linear independence means no cycle. 

Now we generalize the Ehrhart polynomial of $Z_A$. Let $q_1, q_2, \cdots, q_n$ be positive integers. We count the number lattice points in the  $(q_1,\cdots,q_n)$-dilate of $Z_A$, i.e.,
the zonotope $$\widehat{Z_A}:=\{\sum_{i=1}^{n}c_iq_iv_i:0\leq c_i \leq 1\}.$$

\begin{proposition}\label{Ehrhart2}
If $A$ is a totally unimodular matrix, then
$$\#(\widehat{Z_A}\cap\mathbb{Z}^r)=\sum_F\prod_{e\in F}q_e,$$
where $F$ ranges over all the subsets of $\{1,2,\cdots,n\}$ such that $\{v_e:e\in F\}$ is linearly independent. In particular, 
$$\#(\widehat{Z_D}\cap\mathbb{Z}^r)=\sum_{\substack{F\subseteq E\\ F\text{ is acyclic}}}\prod_{e\in F}q_e,$$
\end{proposition}
\begin{proof}
We apply Lemma~\ref{Ehrhart1} to the matrix obtained by multiplying $i$-th column of $A$ by $q_i$ for $i=1,2,\cdots,n$, and take $q=1$.  
\end{proof}

We have a similar counting result for the $(q_1,\cdots,q_n)$-dilate of $S_\sigma$, i.e., $$\widehat{S_\sigma}:=\{(q_1x_1,\cdots,q_nx_n):(x_1,\cdots,x_n)\in S_{\sigma}\}.$$

\begin{proposition}\label{Ehrhart3}
$$\#(\widehat{S_\sigma}\cap\mathbb{Z}^E)=\sum_F\prod_{e\in F}q_e,$$
where $F$ ranges over all the generating forests of the non-empty standard half-open cells $(\text{hoc}_{\sigma,\sigma^*}(T)-\overrightarrow{C^*})\cap [0,1]^E$ in Proposition~\ref{hocprop2}. 
\end{proposition}
\begin{proof}
By Proposition~\ref{hocprop2}, $S_\sigma$ is decomposed into certain standard half-open cells. Note that each of them contains a unique lattice point. So after the $(q_1,\cdots,q_n)$-dilation, each dilated standard half-open cell contains $\prod_{e\in F}q_e$ lattice points, where $F$ is the generating forest of the cell. 
\end{proof}

If the equality $\#(\widehat{Z_D}\cap\mathbb{Z}^r)=\#(\widehat{S_\sigma}\cap\mathbb{Z}^E)$ holds (Proposition~\ref{cc}), then it follows immediately that the map $\varphi_{\sigma,\sigma^*}$ in Proposition~\ref{forestgeometry} is bijective (Theorem~\ref{dd}). Unfortunately, we cannot find a short proof for the equality. To prove it, we need to introduce some definitions and prove some lemmas. 

We call the continuous orientations in $\prod_{e\in E}[0,q_e]$ the \emph{extended} continuous orientations. We restrict the linear map $D$ to the extended continuous orientations. Note that the image is the zonotope $\widehat{Z_D}$. We call this map $\widehat{\psi}$, i.e.,
\begin{eqnarray*}
\widehat{\psi}:\prod_{e\in E}[0,q_e] & \longrightarrow & \widehat{Z_D} \\
\overrightarrow{O} & \mapsto & D\cdot \overrightarrow{O}.
\end{eqnarray*}
Clearly $\widehat{\psi}$ extends $\psi$.

We plan to prove that $\widehat{S_\sigma}$ is a section of $\widehat{\psi}$. To do this, we need to define the \emph{extended $\sigma$-compatible continuous orientations} and prove a similar result (Lemma~\ref{extendedsigmacompatible}) to Proposition~\ref{sigmacompatible2}. The method is the same as in \cite[Prop. 3.2.1]{BBY}.  

Let $\sigma$ be an acyclic cycle signature. An extended continuous orientation $\overrightarrow{O}$ is called \emph{$\sigma$-compatible} if one cannot flow $\overrightarrow{O}$ along any directed cycle in $\sigma$, i.e, $\overrightarrow{O}+\epsilon\sigma(C)\notin \prod_{e\in E}[0,q_e]$ for any $\epsilon >0$ and any cycle $C$. 

Now we prove that the extended $\sigma$-compatible continuous orientations form a section of $\widehat{\psi}$. 
\begin{lemma}\label{extendedsigmacompatible}
For any extended continuous orientation $\overrightarrow{O}$, there exists a unique extended $\sigma$-compatible continuous orientation $\overrightarrow{O^\sigma}$ such that $\overrightarrow{O^\sigma}-\overrightarrow{O}\in\ker_\mathbb{R}(D)$.
\end{lemma}

\begin{proof}
By Lemma~\ref{signature}, there exists $w\in\mathbb{R}^E$ such that $\langle w, \sigma(C)\rangle>0$ for each cycle $C$. Consider the function $P(\overrightarrow{O'}):=\langle w, \overrightarrow{O'}\rangle$ defined on $L:=(\overrightarrow{O}+\ker_\mathbb{R}(D))\cap\prod_{e\in E}[0,q_e]$. If $\overrightarrow{O'}\in L$ is not $\sigma$-compatible, then for some $\epsilon >0$ and some cycle $C$,  $\overrightarrow{O''}:=\overrightarrow{O'}+\epsilon\sigma(C)\in L$ and $P(\overrightarrow{O''})>P(\overrightarrow{O'})$. Because the set $L$ is compact and $P$ is continuous, there exists a maximizer $\overrightarrow{O^\sigma}\in L$ and hence $\overrightarrow{O^\sigma}$ is $\sigma$-compatible. 

It remains to show the uniqueness. Assume by contradiction that there are two such distinct extended continuous orientations $\overrightarrow{O_1}$ and $\overrightarrow{O_2}$. Then $\overrightarrow{O_1}-\overrightarrow{O_2}\in\ker_\mathbb{R}(D)$ and hence by Lemma~\ref{conformal0}, $\overrightarrow{O_1}-\overrightarrow{O_2}=\sum k_i\overrightarrow{C_i}$ where for each edge $e$ of each cycle $C_i$, the sign of $e$ in $\overrightarrow{C_i}$ agrees with the sign of $e$ in $\overrightarrow{O_1}-\overrightarrow{O_2}$. Therefore $\overrightarrow{O_1}-k_1\overrightarrow{C_1}\in D$ and $\overrightarrow{O_2}+k_1\overrightarrow{C_1}\in D$. This contradicts that either $\overrightarrow{C_1}\in\sigma$ or $-\overrightarrow{C_1}\in\sigma$. 
\end{proof}

Then we prove that the set of extended $\sigma$-compatible continuous orientations is exactly $\widehat{S_\sigma}$. The following observation follows immediately by definition. 
\begin{lemma}\label{aa}
Let $\overrightarrow{O}=(x_1,\cdots,x_n)$ be a continuous orientation and $\overrightarrow{O'}=(q_1x_1,\cdots,q_nx_n)$ be an extended continuous orientation. Then $\overrightarrow{O}$ is a $\sigma$-compatible continuous orientation (in the usual sense) if and only if $\overrightarrow{O'}$ is an extended $\sigma$-compatible continuous orientation. 
\end{lemma}

\begin{corollary}\label{bb}
$\widehat{S_\sigma}$ is the set of extended $\sigma$-compatible continuous orientations. 
\end{corollary}
\begin{proof}
Recall that $S_\sigma$ is the set of  $\sigma$-compatible continuous orientations, so this is a direct consequence of Lemma~\ref{aa}. 
\end{proof}

\begin{proposition}\label{cc}
The map $\widehat{\psi}$ restricted to $\widehat{S_\sigma}$ is a bijection to the zonotope $\widehat{Z_D}$. Moreover, $\overrightarrow{O}\in\widehat{S_\sigma}$ is a lattice point if and only if $\widehat{\psi}(\overrightarrow{O})$ is a lattice point. In particular, $$\#(\widehat{Z_D}\cap\mathbb{Z}^r)=\#(\widehat{S_\sigma}\cap\mathbb{Z}^E).$$
\end{proposition}
\begin{proof}
The first part is a direct consequence of Lemma~\ref{extendedsigmacompatible} and Corollary~\ref{bb}. The second part is essentially because $D$ is totally unimodular. Here we give an alternative proof. 

The ``only if'' part is trivial. For the ``if'' part, we consider shifting $\overrightarrow{O}$ into $S_\sigma$. By Proposition~\ref{verticesofS}, $\overrightarrow{O}$ is in the $(q_1,\cdots,q_n)$-dilate of the face $\text{Face}_\sigma(T)$ for some spanning tree $T$. Hence for the edges $e\notin T$, the $e$-th coordinate of $\overrightarrow{O}$ is either $0$ or $q_e$ depending on $\sigma$. So there exists an integer vector $\overrightarrow{v}$ such that the $e$-th coordinate of $\overrightarrow{O}-\overrightarrow{v}$ is either $0$ or $1$ accordingly and the other coordinates are in $[0,1]$, which implies $\overrightarrow{O}-\overrightarrow{v}\in\text{Face}_\sigma(T)$. Because $\widehat{\psi}(\overrightarrow{O}-\overrightarrow{v})$ is a lattice point, by Corollary~\ref{section}, $\overrightarrow{O}-\overrightarrow{v}$ is a lattice point and hence so is $\overrightarrow{O}$.

\end{proof}

Now we are ready to prove the bijectivity of $\varphi_{\sigma,\sigma^*}$.
\begin{theorem}\label{dd}
The map $\varphi_{\sigma,\sigma^*}$ in Proposition~\ref{forestgeometry} is bijective.
\end{theorem}
\begin{proof}
By Proposition~\ref{Ehrhart2}, Proposition~\ref{Ehrhart3}, and Proposition~\ref{cc},
$$\sum_{F_1}\prod_{e\in F_1}q_e=\sum_{F_2}\prod_{e\in F_2}q_e,$$
where $F_1$ ranges over all the spanning forests of $G$ and $F_2$ ranges over all the generating forests of the non-empty standard half-open cells $(\text{hoc}_{\sigma,\sigma^*}(T)-\overrightarrow{C^*})\cap [0,1]^E$ in Proposition~\ref{hocprop2}. Because the two sides are equal for any $(q_1,\cdots,q_n)\in\mathbb{Z}^n$, they are equal as polynomials in $(q_1,\cdots,q_n)$. Hence these generating forests are exactly all the spanning forests of $G$. Therefore the map $\varphi_{\sigma,\sigma^*}$ in Proposition~\ref{forestgeometry} is bijective. 
\end{proof}

The geometric proof of Theorem~\ref{th}(2) is complete. 

\subsection{Duality and Theorem~\ref{th}(3)}

We may get Theorem~\ref{th}(3) by a dual argument. To be precise, for all statements and proofs, we switch cycles and cocycles, $\ker_\mathbb{R}(D)$ and $\im_\mathbb{R}(D^T)$, $\sigma$ and $\sigma^*$, $e\in T$ and $e\notin T$ for a spanning tree $T$, etc.
In particular, for each spanning tree $T$, the face $\text{Face}_{\sigma^*}(T)$ of $[0,1]^E$ is defined to be the set of the continuous orientations where each edge $e\in T$ is oriented according to $\sigma^*(C^*(T,e))$. The tricky part is to replace the matrix $D$ with some matrix $D_*$ and hence get the zonotope $Z_{D_*}:=D_*([0,1]^E)$. The following construction is classic and we refer the readers to \cite[Section 2]{SW} for the technical details. Without loss of generality, we assume the first $r$ columns of $D$ form a basis of the column space of $D$. Hence they form an invertible matrix $P$ and then $P^{-1}D$ is of the form 
\[D_0:=\left[ I_r\: L\right].\]
Note that $\ker_\mathbb{R}(D)=\ker_\mathbb{R}(D_0)$, $\im_\mathbb{R}(D^T)=\im_\mathbb{R}(D_0^T)$, and $D_0$ is totally unimodular. (We can use $D_0$ instead of $D$ in our theory and nothing else needs to be changed.) Now we introduce the ``dual'' matrix. Take \[D_*:=\left[ -L^T \: I_{n-r}\right].\]
It is easy to check $\im_\mathbb{R}(D^T_*)=\ker_\mathbb{R}(D_0)$, $\ker_\mathbb{R}(D_*)=\im_\mathbb{R}(D^T_0)$, and $D_*$ is totally unimodular. Then we can finish the dual argument.

By the dual argument, we obtain a dual geometric construction to the one in Proposition~\ref{hocprop2}, a half-open decomposition of the set $S_{\sigma^*}$ of the $\sigma^*$-compatible continuous orientations. Here we point out that although the maps in Theorem~\ref{th}(2) and (3) agree on the set of $\{(\sigma,\sigma^*)\text{-compatible orientations}\}$, the two geometric constructions do not share the half-open cells corresponding to the trees in general. Indeed,  $\dim(\text{Face}_{\sigma}(T))$ is not equal to $\dim(\text{Face}_{\sigma^*}(T))$ in general and hence the two half-open cells $\text{hoc}_{\sigma,\sigma^*}(T)$ are different. We will see in Section~\ref{tiling} (Proposition~\ref{cor5} and Proposition~\ref{dualcor5}) that they are dual to each other. This is one reason that we cannot combine the two geometric constructions in a naive way and hence derive Theorem~\ref{th}(1).

\section{A geometric interpretation of the main theorem and its combinatorial proof }\label{tiling}

In Section~\ref{geometric}, we have seen that the geometric construction behind Theorem~\ref{th}(2) is the half-open decomposition of the set $S_\sigma$ of the $\sigma$-compatible continuous orientations (Proposition~\ref{forestgeometry}). In this section, we will extend this construction to the cube $[0,1]^E$ by decomposing it into half-open cells (Theorem~\ref{tilingofcube}). This can be viewed as a geometric interpretation of Theorem~\ref{th}(1). 

Our proofs do not rely on the results in Section~\ref{geometric}. We do not make use of zonotopes or translations. Although we need to apply the Ehrhart theory, we still call our proofs combinatorial.

By restricting the half-open decomposition of the cube to $S_\sigma$, we recover the geometric construction in Proposition~\ref{forestgeometry} and its dual construction; see Proposition~\ref{cor5} and Proposition~\ref{dualcor5}. 

\subsection{Half-open decompositions of the cube $[0,1]^E$}
We fix a graph $G$ with the edge set $E$ and a reference orientation.
In this subsection, we describe under which condition a map $$\varphi:\{\text{discrete orientations}\}\longrightarrow\{\text{spanning subgraphs}\}$$ induces a half-open decomposition of the cube $[0,1]^E$ in a canonical way, where $[0,1]^E$ is viewed as the set of continuous orientations of $G$.

Recall the definitions for the \emph{standard half-open cells} given in Section~2.5. For any discrete orientation $\overrightarrow{O}$ and any spanning subgraph $S$, we define the standard half-open cell  $\text{hoc}(\overrightarrow{O},S)$ in the cube $[0,1]^E$ to be the set of the continuous orientations where each edge $e\notin S$ is oriented according to $\overrightarrow{O}$ and each edge $e\in S$ is oriented according to $\overrightarrow{O}$ or bi-oriented. Note that $\text{hoc}(\overrightarrow{O},S)$ contains a unique lattice point $\overrightarrow{O}$ and the generating set is $S$. Conversely, a standard half-open cell in $[0,1]^E$ satisfying these two properties must be $\text{hoc}(\overrightarrow{O},S)$. Hence we identify the standard half-open cells in $[0,1]^E$ with elements in $\{\text{discrete orientations}\}\times\{\text{spanning subgraphs}\}$. We ask which relation, i.e., subset of this Cartesian product, gives a half-open decomposition of $[0,1]^E$. Clearly each discrete orientation should be covered once and only once, so the relation must be a map $\varphi:\{\text{discrete orientations}\}\longrightarrow \{\text{spanning subgraphs}\}$. The target is to study for which maps $\varphi$,  $[0,1]^E=\biguplus_{\overrightarrow{O}\in\{0,1\}^E}\text{hoc}(\overrightarrow{O},\varphi(\overrightarrow{O}))$ holds.

\begin{example}
This example illustrates the half-open decomposition of the cube $[0,1]^E$ induced by a map $\varphi:\{\text{discrete orientations}\}\longrightarrow \{\text{spanning subgraphs}\}$; see Figure~\ref{F6}(a). Here we choose $\varphi$ to be the map $\varphi_{\sigma,\sigma^*}$ in Example~\ref{ex1}; see also Figure~\ref{F5}. We still draw the orientation $\overrightarrow{O}$ and the subgraph $S$ together provided $\varphi_{\sigma,\sigma^*}(\overrightarrow{O})=S$. Now we give this configuration a geometric interpretation: $\text{hoc}(\overrightarrow{O},S)$.
The cube is decomposed into the eight half-open cells $\{\text{hoc}(\overrightarrow{O},\varphi_{\sigma, \sigma^*}(\overrightarrow{O})):\overrightarrow{O}\in\{0,1\}^3\}$. The dimension of each half-open cell is equal to the number of edges in the corresponding subgraph. Here we have $\binom{3}{i}$ half-open cells of dimension $i$ for $i=0,1,2,3$. We will see in Theorem~\ref{tilingofcube} that any $\varphi_{\sigma,\sigma^*}$ in Theorem~\ref{th} induces a half-open decomposition. 

\begin{figure}[h]
            \centering
            \includegraphics[scale=1]{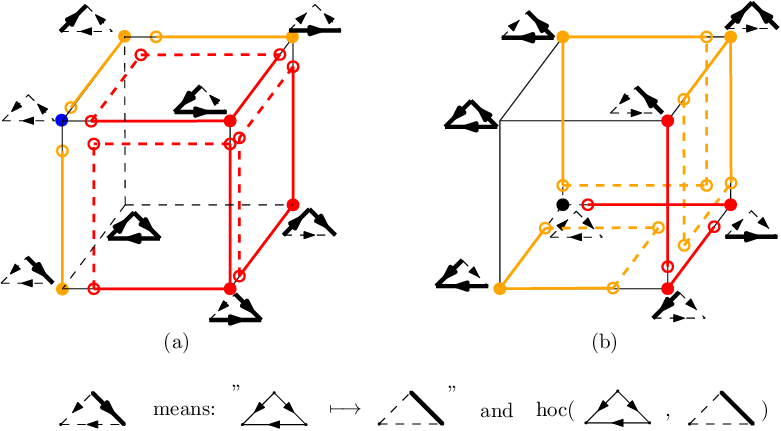}
            \caption{(a) This shows $[0,1]^E=\biguplus_{\protect\overrightarrow{O}\in\{0,1\}^E}\text{hoc}(\protect\overrightarrow{O},\varphi_{\sigma,\sigma^*}(\protect\overrightarrow{O}))$, where the bijection $\varphi_{\sigma,\sigma^*}$ is as in Figure~\ref{F5}. The unique 3-dimensional half-open cell is not  drawn to keep the figure clean. (b) This is dual tiling of (a); see the definition of the dual tiling after Corollary~\ref{cordual}. The unique 3-dimensional half-open cell is not drawn.}
            \label{F6}
\end{figure}

\end{example}

The main result of the subsection is as follows. 

\begin{proposition}\label{abcd}
For a map $\varphi:\{\text{discrete orientations}\}\longrightarrow \{\text{spanning subgraphs}\}$, the following are equivalent. 

(\ref{abcd}a) $[0,1]^E=\biguplus_{\overrightarrow{O}\in\{0,1\}^E}\text{hoc}(\overrightarrow{O},\varphi(\overrightarrow{O}))$.

(\ref{abcd}b) $\varphi$ is bijective and for any two distinct discrete orientations $\overrightarrow{O_1}$ and 
$\overrightarrow{O_2}$, $\text{hoc}(\overrightarrow{O_1},\varphi(\overrightarrow{O_1}))\cap\text{hoc}(\overrightarrow{O_2},\varphi(\overrightarrow{O_2}))=\emptyset$.

(\ref{abcd}c) For any two distinct discrete orientations $\overrightarrow{O_1}$ and 
$\overrightarrow{O_2}$, there exists an edge $e$ such that $\overrightarrow{O_1}(e)\neq\overrightarrow{O_2}(e)$ and $e\in \varphi(\overrightarrow{O_1})\bigtriangleup \varphi(\overrightarrow{O_2})$.

(\ref{abcd}d) $\varphi$ is \emph{locally bijective} (defined as follows similarly to Proposition~\ref{local}(1)).

\end{proposition}

\begin{definition}
A map $\varphi:\{\text{discrete orientations}\}\longrightarrow \{\text{spanning subgraphs}\}$ is called \emph{locally bijective} if for any partial orientation $\overrightarrow{O_P}$ of $G$, where $P$ is the underlying edges of $\overrightarrow{O_P}$, the map 

\begin{eqnarray*}
\varphi_{\overrightarrow{O_P}}:\{\text{orientations of }E\backslash P\} & \longrightarrow & \{\text{subsets of }E\backslash P\} \\
\overrightarrow{O_{E\backslash P}} & \mapsto & \varphi(\overrightarrow{O_{E\backslash P}}\cup\overrightarrow{O_P})\backslash P
\end{eqnarray*}
is bijective.
\end{definition}

Remark that we can also add the restriction in Proposition~\ref{local}(2) as statement (\ref{abcd}e), but it does not help proving the equivalences among (\ref{abcd}a), (\ref{abcd}b), (\ref{abcd}c), and (\ref{abcd}d).

First we prove (\ref{abcd}a) $\Leftrightarrow$ (\ref{abcd}b) by using the generalized Ehrhart polynomial.  

\begin{lemma}\label{ab}
(\ref{abcd}a) $\Leftrightarrow$ (\ref{abcd}b).
\end{lemma}

\begin{proof}
We need to prove that if the second part of (\ref{abcd}b) holds, then $\varphi$ is bijective if and only if (\ref{abcd}a) holds. For any positive integers $q_1, \cdots, q_n$, where $n=|E|$, we consider the $(q_1,\cdots,q_n)$-dilation:

\begin{eqnarray*}
f: \mathbb{R}^E & \longrightarrow & \mathbb{R}^E \\
(x_1, \cdots, x_n) & \mapsto & (q_1x_1, \cdots, q_nx_n).
\end{eqnarray*}

Dilating the cube, we get $$\#( f([0,1]^E)\cap\mathbb{Z}^E)=\prod_{e\in E}(1+q_e).$$ 

Dilating the half-open cells, we get $$\#( f(\biguplus_{\overrightarrow{O}\in\{0,1\}^E}\text{hoc}(\overrightarrow{O},\varphi(\overrightarrow{O})))\cap\mathbb{Z}^E)=\sum_{\overrightarrow{O}\in\{0,1\}^E}\prod_{e\in\varphi(\overrightarrow{O})}q_e.$$ 

Consider the equation
\begin{equation}\label{ast}
\#( f([0,1]^E)\cap\mathbb{Z}^E)=\#( f(\biguplus_{\overrightarrow{O}\in\{0,1\}^E}\text{hoc}(\overrightarrow{O},\varphi(\overrightarrow{O})))\cap\mathbb{Z}^E).
\end{equation}

By the previous two formulas, one gets that Equation~\ref{ast} holds for any $(q_1,\cdots,q_n)$-dilation $f$ $\Leftrightarrow$ $\prod\limits_{e\in E}(1+q_e)=\sum\limits_{\overrightarrow{O}\in\{0,1\}^E}\prod\limits_{e\in\varphi(\overrightarrow{O})}q_e$ as polynomials in $q_e$'s $\Leftrightarrow$ $\{\varphi(\overrightarrow{O}):\overrightarrow{O}\in\{0,1\}^E\}=\{\text{spanning subgraphs}\}$ $\Leftrightarrow$ $\varphi$ is bijective.

It remains to show that (\ref{abcd}a) $\Leftrightarrow$ Equation~\ref{ast} holds for any $(q_1,\cdots,q_n)$-dilation $f$. The ``$\Rightarrow$'' part is trivial. Now we begin to prove the other direction. Clearly $\text{hoc}(\overrightarrow{O},\varphi(\overrightarrow{O}))\subseteq [0,1]^E$. Assume by contradiction that some point in the cube is not covered by the half-open cells. Then one can disturb the point such that it is a rational point and is still not covered. By a suitable choice of $f$, this rational point becomes integral after the dilation, which contradicts the cardinality equality. Therefore (\ref{abcd}a) $\Leftrightarrow$ (\ref{abcd}b).
\end{proof}

\begin{lemma}\label{acd}
(\ref{abcd}c) $\Rightarrow$ (\ref{abcd}b) $\Rightarrow$ (\ref{abcd}a) $\Rightarrow$ (\ref{abcd}d) $\Rightarrow$ (\ref{abcd}c).
\end{lemma}
\begin{proof}

It is trivial that (\ref{abcd}c) implies the injectivity of  $\varphi$ and hence the bijectivity. It is straightforward to check that (\ref{abcd}c) implies the second part of (\ref{abcd}b). So (\ref{abcd}c) $\Rightarrow$ (\ref{abcd}b). 

(\ref{abcd}b) $\Rightarrow$ (\ref{abcd}a) is proved in Lemma~\ref{ab}.

For (\ref{abcd}a) $\Rightarrow$ (\ref{abcd}d), we consider the restriction of the tiling. Observe that for a partial orientation  $\overrightarrow{O_P}$, the set of the continuous orientations extending $\overrightarrow{O_P}$ is a closed face $F$ of $[0,1]^E$. We restrict the decomposition  $[0,1]^E=\biguplus_{\overrightarrow{O}\in\{0,1\}^E}\text{hoc}(\overrightarrow{O},\varphi(\overrightarrow{O}))$ to the closed face $F$.  Then we get $$F=[0,1]^{E\backslash P}\times\overrightarrow{O_P}|_P =\biguplus_{\overrightarrow{O}\in\{0,1\}^E \text{ extends } \overrightarrow{O_P}}\text{hoc}(\overrightarrow{O},\varphi(\overrightarrow{O}))\cap F.$$
It is straightforward to check that for any $\overrightarrow{O}\in\{0,1\}^E$ that extends  $\overrightarrow{O_P}$, we have
$$\text{hoc}(\overrightarrow{O},\varphi(\overrightarrow{O}))\cap F=\text{hoc}(\overrightarrow{O}|_{E\backslash P},\varphi(\overrightarrow{O})\backslash P)\times\overrightarrow{O_P}|_P=\text{hoc}(\overrightarrow{O}|_{E\backslash P},\varphi_{\overrightarrow{O_P}}(\overrightarrow{O}|_{E\backslash P}))\times\overrightarrow{O_P}|_P,$$
where the later two ``hoc'' are with respect to the cube $[0,1]^{E\backslash P}$. Hence 
$$[0,1]^{E\backslash P}=\biguplus_{\overrightarrow{O_{E\backslash P}}\in\{0,1\}^{E\backslash P}}\text{hoc}(\overrightarrow{O_{E\backslash P}},\varphi_{\overrightarrow{O_P}}(\overrightarrow{O_{E\backslash P}})).$$
By applying Lemma~\ref{ab} to $\varphi_{\overrightarrow{O_P}}$, we get $\varphi_{\overrightarrow{O_P}}$ is bijective and hence (\ref{abcd}d) holds. 

It remains to show (\ref{abcd}d) $\Rightarrow$ (\ref{abcd}c). Denote $E_{\rightrightarrows}=\{e:\overrightarrow{O_1}(e)=\overrightarrow{O_2}(e)\}$ and $E_{\rightleftarrows}=\{e:\overrightarrow{O_1}(e)\neq\overrightarrow{O_2}(e)\}$($\neq\emptyset$). By (\ref{abcd}d), 
$\varphi_{\overrightarrow{O_{E_{\rightrightarrows}}}}$ is bijective. Hence $\varphi_{\overrightarrow{O_{E_{\rightrightarrows}}}}(\overrightarrow{O_1}|_{E_{\rightleftarrows}})\neq \varphi_{\overrightarrow{O_{E_{\rightrightarrows}}}}(\overrightarrow{O_2}|_{E_{\rightleftarrows}})$. By the definition of $\varphi_{\overrightarrow{O_{E_{\rightrightarrows}}}}$, this means  $\varphi(\overrightarrow{O_1})\backslash E_{\rightrightarrows}\neq\varphi(\overrightarrow{O_2})\backslash E_{\rightrightarrows}$. So there exists an edge $e\in E_{\rightleftarrows}$ such that $e\in \varphi(\overrightarrow{O_1})\bigtriangleup \varphi(\overrightarrow{O_2})$. So (\ref{abcd}c) holds. 
\end{proof}

By Lemma~\ref{ab} and Lemma~\ref{acd}, Proposition~\ref{abcd} holds.

When we define $\text{hoc}(\overrightarrow{O},S)$, 
it is somewhat arbitrary that we bi-orient the edges in $S$ rather than the edges in $E\backslash S$. Actually, we could also define it in the other way, or equivalently, we may replace $\varphi$ with $\varphi^*$ defined by $\overrightarrow{O}\mapsto E\backslash\varphi(\overrightarrow{O})$. 

\begin{corollary}\label{cordual}
If (\ref{abcd}a) holds for $\varphi$, then (\ref{abcd}a) also holds for $\varphi^*$, where the map $\varphi^*$ is defined by $\overrightarrow{O}\mapsto E\backslash\varphi(\overrightarrow{O})$. 
\end{corollary}
\begin{proof}
By Proposition~\ref{abcd}, (\ref{abcd}a) $\Leftrightarrow$ (\ref{abcd}c). Clearly (\ref{abcd}c) holds for $\varphi$ if and only if (\ref{abcd}c) holds for $\varphi^*$, so (\ref{abcd}a) holds for $\varphi^*$.
\end{proof}

We call the tiling $[0,1]^E=\biguplus_{\overrightarrow{O}\in\{0,1\}^E}\text{hoc}(\overrightarrow{O},\varphi^*(\overrightarrow{O}))$ the \emph{dual tiling} of $[0,1]^E=\biguplus_{\overrightarrow{O}\in\{0,1\}^E}\text{hoc}(\overrightarrow{O},\varphi(\overrightarrow{O}))$. Clearly $(\varphi^*)^*=\varphi$. See Figure~\ref{F6}(b) for an example of the dual tiling. 

\begin{Rem}
Although we develop the theory in this subsection for a graph $G$, the structure of $G$ does not play a role in either the results or the proofs. What we need is just a finite set $E$, whose elements are viewed as edges.
\end{Rem}

\subsection{A geometric interpretation of Theorem~\ref{th}}

By Proposition~\ref{propsubset}, $\varphi_{\sigma,\sigma^*}$ satisfies (\ref{abcd}c). By the equivalence of (\ref{abcd}a) and (\ref{abcd}c) in Proposition~\ref{abcd}, we have the following result. 

\begin{theorem}\label{tilingofcube}
$[0,1]^E=\biguplus_{\overrightarrow{O}\in\{0,1\}^E}\text{hoc}(\overrightarrow{O},\varphi_{\sigma,\sigma^*}(\overrightarrow{O}))$.
\end{theorem}

Recall that if we choose the map $\varphi_{\sigma,\sigma^*}$ in Example~\ref{ex1}, then we get the decomposition in Figure~\ref{F6}(a). 

The tiling of the cube $[0,1]^E$ in Theorem~\ref{tilingofcube} can be viewed as a geometric construction behind the bijection $\varphi_{\sigma,\sigma^*}$ in Theorem~\ref{th}(1) for the following two reasons. First, the bijection $\varphi_{\sigma,\sigma^*}$ induces the tiling, and conversely, from the tiling, we may get the map $\varphi_{\sigma,\sigma^*}$ back by sending an orientation $\overrightarrow{O}$ to the generating set of the unique standard half-open cell containing $\overrightarrow{O}$. Second, the tiling extends the tiling of $S_\sigma$ in Proposition~\ref{hocprop2}, which is viewed as the geometric construction behind the bijection in Theorem~\ref{th}(2). Moreover, the dual tiling extends the tiling of $S_{\sigma^*}$ behind Theorem~\ref{th}(3). We now begin to show the second point.

We consider the restriction of the tiling in Theorem~\ref{tilingofcube} to $S_\sigma$, the set of $\sigma$-compatible continuous orientations.
We will prove the following decomposition of $S_\sigma$ and its dual version, Proposition~\ref{dualcor5}.
\begin{proposition}\label{cor5}
$S_\sigma=\biguplus_{\overrightarrow{O}\in\{0,1\}^E\text{ is }\sigma\text{-compatible}}\text{hoc}(\overrightarrow{O},\varphi_{\sigma,\sigma^*}(\overrightarrow{O}))$.
\end{proposition}
\begin{Rem}\label{rmk}
Comparing Proposition~\ref{forestgeometry} and Proposition~\ref{cor5}, we have the same map $\varphi_{\sigma,\sigma^*}$. Moreover, the half-open cells in Proposition~\ref{hocprop2} and  Proposition~\ref{cor5} induce the map $\varphi_{\sigma,\sigma^*}$ in the same way, which sends the unique lattice point contained in a half-open cell to  the generating set. So the half-open decomposition in Proposition~\ref{hocprop2} is the same as the one in Proposition~\ref{cor5}. Hence the decomposition in Theorem~\ref{tilingofcube} extends the one in Proposition~\ref{hocprop2}. For example, the tiling in Figure~\ref{F6} extends the tiling in Figure~\ref{F3}(b). By Proposition~\ref{forestgeometry}, we already know Proposition~\ref{cor5} holds. However, the proof we will give in this section is combinatorial and independent of the proofs in Section~\ref{geometric}, which makes our combinatorial approach self-contained. 
\end{Rem}

To prove Proposition~\ref{cor5}, we need the following lemma, which is part of Proposition~\ref{verticesofS} and is proved implicitly in \cite{BBY}. For readers' convenience, we give a proof here. Recall that $\text{Face}_\sigma(T)$ is a closed face of $[0,1]^E$ consisting of the continuous orientations where each edge $e\notin T$ is oriented according to $\sigma(C(T,e))$. 

\begin{lemma}\label{facelemma}
$\text{Face}_\sigma(T)\subseteq S_\sigma$. 
\end{lemma}
\begin{proof}
For any continuous orientation $\overrightarrow{O}\in\text{Face}_\sigma(T)$, we need to show it is $\sigma$-compatible, i.e., for any $\epsilon >0$ and any cycle $C$, $\overrightarrow{O}+\epsilon\sigma(C)\notin [0,1]^E$. By Lemma~\ref{fundamental}, $\sigma(C)=\sum\limits_{e\notin T,\overrightarrow{e}\in\sigma (C)}C(T,\overrightarrow{e})$. Because $\sigma(C)$ is $\sigma$-compatible, at least one of the directed cycles in the sum is $\sigma$-compatible, say $C(T,\overrightarrow{e_0})$. Note that $\overrightarrow{e_0}$ is an arc in $\sigma(C)$ and by the definition of $\text{Face}_\sigma(T)$,  $\overrightarrow{e_0}$ is an arc in $\overrightarrow{O}$. So $\overrightarrow{O}+\epsilon\sigma(C)\notin [0,1]^E$. 
\end{proof}

We also need the following notation and lemma. Let $\overrightarrow{O}$ be a continuous orientation. Denote some of the discretely oriented edges in $\overrightarrow{O}$ by a partial orientation $\overrightarrow{P}$. Then we say $\overrightarrow{O}$ contains $\overrightarrow{P}$ and denote by $_{\overrightarrow{P}}\overrightarrow{O}$ or $_{P}\overrightarrow{O}$ the continuous orientation obtained by reversing the arcs $\overrightarrow{P}$ in $\overrightarrow{O}$.

\begin{lemma}\label{lemmacc}
Let $\overrightarrow{O}$ be a $\sigma$-compatible continuous orientation. If $\overrightarrow{O}$ contains a directed cocycle $\overrightarrow{C^*}$, then $_{C^*}\overrightarrow{O}$ is also $\sigma$-compatible. 
\end{lemma}

\begin{proof}
Assume by contradiction that there exists some $\epsilon>0$ and cycle $C$ such that $_{C^*}\overrightarrow{O}+\epsilon\sigma(C)\in [0,1]^E$. 

Note that $C\cap C^*$ must be empty. Indeed, if $C\cap C^*\neq\emptyset$, because $\langle\sigma(C),\overrightarrow{C^*}\rangle=0$, one arc in $\sigma(C)$ will be equal to one arc in $\overrightarrow{C^*}$ and some other arc in $\sigma(C)$ will be opposite to an arc in $\overrightarrow{C^*}$. Hence $_{C^*}\overrightarrow{O}+\epsilon\sigma(C)\notin [0,1]^E$, which leads to a contradiction. 

Because $C\cap C^*=\emptyset$ and $_{C^*}\overrightarrow{O}+\epsilon\sigma(C)\in [0,1]^E$, $\overrightarrow{O}+\epsilon\sigma(C)\in [0,1]^E$, which contradicts that $\overrightarrow{O}$ is  $\sigma$-compatible. 
\end{proof}

Now we show which half-open cells in Theorem~\ref{tilingofcube} consist of $\sigma$-compatible continuous orientations. 

\begin{lemma}\label{lemmadd}
Let $\overrightarrow{O}$ be a discrete orientation. Then a continuous orientation in $\text{hoc}(\overrightarrow{O},\varphi_{\sigma,\sigma^*}(\overrightarrow{O}))$ is $\sigma$-compatible if and only if $\overrightarrow{O}$ is $\sigma$-compatible. 
\end{lemma}
\begin{proof}
Adopt the notations in Lemma~\ref{l1}. 

If $\overrightarrow{O}$ is not $\sigma$-compatible, then by Lemma~\ref{l1}, $I\neq\emptyset$, which means $\overrightarrow{O}$ contains a directed cycle $-\overrightarrow{C_i}\notin\sigma$ for some $i\in I$. Note that $\varphi_{\sigma,\sigma^*}(\overrightarrow{O})$ contains the cycle $C_i$. Hence for any continuous orientation $\overrightarrow{O'}\in\text{hoc}(\overrightarrow{O},\varphi_{\sigma,\sigma^*}(\overrightarrow{O}))$, there exist a sufficiently small $\epsilon$ such that $\overrightarrow{O'}+\epsilon\overrightarrow{C_i}\in [0,1]^E$. Therefore $\overrightarrow{O'}$ is not $\sigma$-compatible. 

If $\overrightarrow{O}$ is $\sigma$-compatible, then by Lemma~\ref{l1}, $I=\emptyset$. Let $T$ be the spanning tree $\text{BBY}_{\sigma,\sigma^*}^{-1}(\overrightarrow{O^{cp}})$. Note that $\overrightarrow{O}=_{\biguplus_{j\in J}C^*_j}\overrightarrow{O^{cp}}$ and $\varphi_{\sigma,\sigma^*}(\overrightarrow{O})=T\backslash\biguplus_{j\in J}C^*_j$. So any continuous orientation $\overrightarrow{O'}\in\text{hoc}(\overrightarrow{O},\varphi_{\sigma,\sigma^*}(\overrightarrow{O}))$ contains directed cocycle $-\overrightarrow{C^*_j}$'s and $_{\biguplus_{j\in J}C^*_j}\overrightarrow{O'}\in\text{hoc}(\overrightarrow{O^{cp}},T\backslash\biguplus_{j\in J}C^*_j)\subseteq\text{hoc}(\overrightarrow{O^{cp}},T)\subseteq\text{Face}_\sigma(T)$. By Lemma~\ref{facelemma}, $_{\biguplus_{j\in J}C^*_j}\overrightarrow{O'}$ is $\sigma$-compatible. By Lemma~\ref{lemmacc}, $\overrightarrow{O'}$ is also $\sigma$-compatible. 
\end{proof}

By Lemma~\ref{lemmadd} and Theorem~\ref{tilingofcube}, we get Proposition~\ref{cor5}.

Now we consider the dual theory. We only list the results because the proof is similar. 

Due to Corollary~\ref{cordual}, we have the following the dual decomposition. 
\begin{proposition}
$[0,1]^E=\biguplus_{\overrightarrow{O}\in\{0,1\}^E}\text{hoc}(\overrightarrow{O},\varphi^*_{\sigma,\sigma^*}(\overrightarrow{O}))$.
\end{proposition}

Let $S_{\sigma^*}$ be the set of $\sigma^*$-compatible continuous orientations. The face $\text{Face}_{\sigma^*}(T)$ of $[0,1]^E$ is defined to be the set of the continuous orientations where each edge $e\in T$ is oriented according to $\sigma^*(C^*(T,e))$. Then we have the following dual lemmas.

\begin{lemma}
(1) $\text{Face}_{\sigma^*}(T)\subseteq S_{\sigma^*}$. 

(2) Let $\overrightarrow{O}$ be a $\sigma^*$-compatible continuous orientation. If $\overrightarrow{O}$ contains a directed cycle $\overrightarrow{C}$, then $_{C}\overrightarrow{O}$ is also $\sigma^*$-compatible. 

(3) Let $\overrightarrow{O}$ be a discrete orientation. Then a continuous orientation in $\text{hoc}(\overrightarrow{O},\varphi^*_{\sigma,\sigma^*}(\overrightarrow{O}))$ is $\sigma^*$-compatible if and only if $\overrightarrow{O}$ is $\sigma^*$-compatible. 
\end{lemma}

Finally, we get the decomposition of $S_{\sigma^*}$. 
\begin{proposition}\label{dualcor5}
$S_{\sigma^*}=\biguplus_{\overrightarrow{O}\in\{0,1\}^E\text{ is }\sigma^*\text{-compatible}}\text{hoc}(\overrightarrow{O},\varphi^*_{\sigma,\sigma^*}(\overrightarrow{O}))$.
\end{proposition}

\begin{figure}[h]
            \centering
            \includegraphics[scale=1]{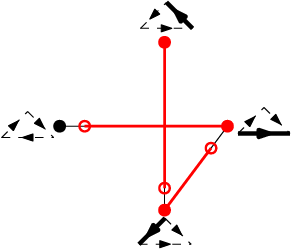}
            \caption{This shows the half-open decomposition of $S_{\sigma^*}$ in Example~\ref{ex}, which induces the bijection $\varphi_{\sigma,\sigma^*}$ from $\sigma$-compatible orientations to connected spanning subgraphs.}
            \label{F7}
\end{figure}

\begin{example}\label{ex}
We still consider the bijection $\varphi_{\sigma,\sigma^*}$ in Example~\ref{ex1}. Note that for a spanning tree $T$, the closed face $\text{Face}_{\sigma^*}(T)$ and the half-open cell $\text{hoc}(\text{BBY}_{\sigma, \sigma^*}(T),E\backslash T)$ ($=\text{hoc}(\overrightarrow{O^{cp}},\varphi^*_{\sigma,\sigma^*}(\overrightarrow{O^{cp}}))$, where $\overrightarrow{O^{cp}}=\text{BBY}_{\sigma, \sigma^*}(T)$) are of dimension $1$. By restricting the tiling in Figure~\ref{F6}(b) to $S_{\sigma^*}$, we get Figure~\ref{F7}, which consists of three half-open line segments and one point. 
\end{example}

\section{Generalization to Regular Matroids}\label{regular}
In this section, we first introduce the definition of regular matroids and related objects; see also \cite{BBY}. Then we explain why the results on graphs in Section~\ref{geometric} and Section~\ref{combinatorial} can be generalized to the regular matroids. At last, we explain the obstacles to generalizing the results to \emph{realizable matroids} over $\mathbb{R}$.  

We assume that the reader is familiar with the basic theory of matroids; some standard references include \cite{O}. Recall that a matrix is called \emph{totally unimodular} if every square submatrix has determinant $0$, $1$, or $-1$. A \emph{matroid} is called \emph{regular} if it can be represented by a totally unimodular matrix over $\mathbb{R}$. Let $A$ be an $r\times n$ matrix representing a regular matroid $M$ with the ground set $E$. Without loss of generality, we may assume $\text{rank}(A)=r$. 

Let $C$ be a circuit. By definition, $C$ is the support of a support-minimal nonzero element in $\ker_\mathbb{R}(A)$. By \cite[Lemma 6]{SW}\footnote{In \cite{SW}, Lemma~6 and Lemma~7 are proved with respect to $\ker_\mathbb{R}(A)\cap\mathbb{Z}^n$ instead of $\ker_\mathbb{R}(A)$. However, the proofs work for $\ker_\mathbb{R}(A)$.}, all the elements in $\ker_\mathbb{R}(A)$ with support $C$, together with the zero vector, form a one-dimensional subspace of $\ker_\mathbb{R}(A)$. By \cite[Lemma 7]{SW}, the generator of this subspace can be chosen to be a $\{0, \pm 1\}$-vector. Clearly there are exactly two such generators, denoted by $\pm\overrightarrow{C}$. So for any circuit $C$, there are exactly two $\{0, \pm 1\}$-vectors in $\ker_\mathbb{R}(A)$ with support $C$. We call them \emph{signed circuits} of $M$. Note that if we take $A$ to be the matrix $D$ associated to a graph $G$ (defined in Section~\ref{combinatorial}), then these signed circuits are exactly the directed cycles. 

By \cite[Lemma 10]{SW}, the notion of signed circuit is intrinsic to $M$, independent of the choice of $A$, and hence we are safe to work with one matrix $A$. To be precise, let $A'$ be another $r\times n$ totally unimodular matrix representing $M$. Without loss of generality, we may assume the $i$-th columns of $A$ and $A'$ correspond to the same element in $E$ for $i=1, \cdots, n$. Then we have $A'=FAP$, where $F$ is an integer matrix whose determinant is $\pm 1$ and $P$ is a diagonal matrix whose entries on the main diagonal are $\pm 1$; see \cite[Section 2.2]{SW}. Hence for any vector $\overrightarrow{C}$, $P\cdot\overrightarrow{C}\in\ker_\mathbb{R}(A)$ if and only if $\overrightarrow{C}\in\ker_\mathbb{R}(A')$. So the signed circuits with respect to $A$ and the ones with respect to $A'$ differ merely by a reorientation $P$ of $M$. As in the case of graphs, when we choose the matrix $A$, it contains the information of a reference orientation for $M$. The choice of the reference orientation does not affect our results essentially. We remark that these signed circuits make $M$ an oriented matroid; see \cite[Chapter 1.2]{BVSWZ}. Moreover, the oriented matroid structure on a regular matroid is unique up to reorientation; see \cite[Corollary 7.9.4]{BVSWZ}. 

Similarly, for any cocircuit $C^*$, there are exactly two $\{0, \pm 1\}$-vectors in $\im_\mathbb{R}(A)$ with support $C^*$. We call them \emph{signed cocircuits} of $M$. All the arguments above for circuits also work for cocircuits due to the duality; 
see \cite[Section 2]{SW} or Subsection~\ref{subsectionduality}. 

One can generalize the notion of fundamental cycles and cocycles, acyclic cycle signature and cocycle signature, cycle reversals, and cocycle reversals in a straightforward way from graphs to regular matroids, while one needs to replace ``cycle'' with ``circuit'' and ``cocycle'' with ``cocircuit'' in the names. For details, see \cite{BBY}. 

Similarly to the case of graphs, the vector space $\mathbb{R}^E$ has an orthogonal decomposition $\ker_\mathbb{R}(A)\oplus\im_\mathbb{R}(A^T)$ with respect to the standard inner product $\langle\cdot,\cdot\rangle$. Hence the signed circuit and the signed cocircuit are orthogonal to each other. For any basis $B$, it is easy to check that the signed fundamental circuits (resp. signed fundamental cocircuits) form a basis of $\ker_\mathbb{R}(A)$ (resp. $\im_\mathbb{R}(A^T)$), and an integral basis of $\ker_\mathbb{Z}(A)=\ker_\mathbb{R}(A)\cap \mathbb{Z}^E$ (resp. $\im_\mathbb{Z}(A^T)=\im_\mathbb{R}(A^T)\cap \mathbb{Z}^E$). 

Now we claim that we may generalize all the results in the previous sections to the regular matroids. Indeed, the results that we cite from other references are proved in the setting of regular matroids in these references, most of which are in \cite{BBY}. For the proofs in this paper, one can simply replace ``cycle'' with ``circuit'',  ``cocycle'' with ``cocircuit'', ``spanning tree'' with ``basis'', etc. 

In particular, Proposition~\ref{tree} and Theorem~\ref{th} hold for the regular matroids.

\begin{proposition}(\cite{BBY},Theorem 1.3.1(1))
Fix acyclic signatures $\sigma$ and $\sigma^*$ of a regular matroid $M$.  
Then the map $$\text{BBY}_{\sigma,\sigma^*}:\{\text{bases of } M\}\longrightarrow\{(\sigma,\sigma^*)\text{-compatible orientations}\}$$ is a bijection, where $\text{BBY}_{\sigma,\sigma^*}$ sends $B$ to the orientation of $M$ in which we orient each $e\notin B$ according to its orientation in $\sigma(C(B,e))$ and each $e\in B$ according to its orientation in $\sigma^*(C^*(B,e))$.
\end{proposition}

\begin{theorem}

Fix acyclic signatures $\sigma$ and $\sigma^*$ of a regular matroid $M$ with ground set $E$.  

(1) The map 

\begin{eqnarray*}
\varphi_{\sigma,\sigma^*}:\{\text{discrete orientations}\} & \longrightarrow & \{\text{subsets of } E\} \\
\overrightarrow{O} & \mapsto & (\text{BBY}_{\sigma,\sigma^*}^{-1}(\overrightarrow{O^{cp}})\cup \biguplus_{i\in I}C_i)\backslash \biguplus_{j\in J}C_j^*
\end{eqnarray*}
is a bijection, where $\overrightarrow{O}$ is an orientation obtained by reversing disjoint signed circuits $\{\overrightarrow{C_i}\}_{i\in I}$ and signed cocircuits $\{\overrightarrow{C_j^*}\}_{j\in J}$ in a ($\sigma,\sigma^*$)-compatible orientation $\overrightarrow{O^{cp}}$.

(2) The map $\varphi_{\sigma,\sigma^*}$ specializes to the bijection
\begin{eqnarray*}
\varphi_{\sigma,\sigma^*}: \{\sigma\text{-compatible orientations}\} & \longrightarrow & \{\text{independent sets}\} \\
\overrightarrow{O} & \mapsto & \text{BBY}_{\sigma,\sigma^*}^{-1}(\overrightarrow{O^{cp}})\backslash \biguplus_{j\in J}C_j^*.
\end{eqnarray*}

(3) The map $\varphi_{\sigma,\sigma^*}$ specializes to the bijection
\begin{eqnarray*}
\varphi_{\sigma,\sigma^*}:\{\sigma^*\text{-compatible orientations}\} & \longrightarrow & \{\text{spanning sets}\} \\
\overrightarrow{O} & \mapsto & \text{BBY}_{\sigma,\sigma^*}^{-1}(\overrightarrow{O^{cp}})\cup \biguplus_{i\in I}C_i.
\end{eqnarray*}
\end{theorem}

Note that Proposition~\ref{tree} is proven for \emph{realizable matroids} over $\mathbb{R}$ in \cite{BBY} and, more generally, for \emph{oriented matroids} in \cite{BSY}. It is  natural to ask whether Theorem~\ref{th} also holds for these types of matroids. The answer is negative (unless some notions are modified significantly). A direct reason is that one circuit-cocircuit reversal class of a realizable matroid possibly contains more than one $(\sigma,\sigma^*)$-compatible orientations. To be precise, we define the circuit reversals and the cocircuit reversals for the realizable matroids by analogy with the terminology for graphs; see also \cite{G2} for the definitions. We assume the $(\sigma,\sigma^*)$-compatible orientations for the realizable matroids are defined so that they are in bijection with the bases; see \cite{BSY} for such a definition. We consider the uniform oriented matroid $U_{2,4}$, which is realizable. It is easy to check that $U_{2,4}$ has only two circuit-cocircuit reversal classes while it has six bases (cf. \cite[Prop. 3.2]{G2}). Hence for an orientation $\overrightarrow{O}$, there may be more than one $(\sigma,\sigma^*)$-compatible orientations $\overrightarrow{O^{cp}}$ in the same circuit-cocircuit reversal class, which makes the formulation of Theorem~\ref{th} problematic.

\begin{Rem}
We also mention a related result here. In \cite[Theorem 3.3]{GY}, Gioan and Yuen proved that, for an oriented matroid, the number of circuit-cocircuit reversal classes equals the number of bases if and only if the matroid is regular. 
\end{Rem}

We would like to share another way to see the obstacles to the generalization. 
From a geometric point of view, the zonotope generated by a set $Z$ of vectors tiles the space if and only if the matroid represented by $Z$ is regular (see \cite[Theorem 1]{DG}). Hence our geometric construction as in Figure~\ref{F3} cannot work beyond the regular matroids.

\section*{Acknowledgement}

Many thanks to Olivier Bernardi for orienting the author towards the study of this project and countless helpful discussions, as well as detailed advice on the writing of this paper. As a graduate student, the author wants to thank the Department of Mathematics at Brandeis University for admitting the author to the Ph.D. program and creating a flexible study environment. The author also wants to thank Chi Ho Yuen for helpful discussions. Thanks to the anonymous referees for the helpful feedback.

\end{document}